\documentclass[moor]{informs1}              
\usepackage{natbib}
 \NatBibNumeric
 \def\bibfont{\small}%
 \bibpunct[, ]{[}{]}{,}{n}{}{,}%

\usepackage[colorlinks=true,breaklinks=true,bookmarks=true,urlcolor=blue,
     citecolor=blue,linkcolor=blue,bookmarksopen=false,draft=false]{hyperref}

\def\EMAIL#1{\href{mailto:#1}{#1}}

\TheoremsNumberedThrough     

\EquationsNumberedThrough    

\usepackage{mathrsfs}
\usepackage{graphics}
\usepackage{enumerate}

\usepackage{tikz} 
\usetikzlibrary{automata, arrows}

\newtheorem{Assumption}{Assumption}[section]
\newcommand{\verteq}{\rotatebox{90}{$\,=$}}

\def\B{\mathfrak{B}}

\def\BB{\mathbb}
\def\F{\mathscr{F}}
\def\P{\mathscr{P}}
\def\R{\mathbb{R}}

\begin{document}


\RUNAUTHOR{Feinberg, Mandava, and Shiryaev}

\RUNTITLE{Sufficiency of Markov policies for CTJMDPs}

\TITLE{Sufficiency of Markov Policies for Continuous-Time Jump Markov Decision Processes}

\ARTICLEAUTHORS{%
\AUTHOR{Eugene A. Feinberg}
\AFF{Stony Brook University, Stony Brook, New York,  11794, USA, \EMAIL{eugene.feinberg@stonybrook.edu}}
\AUTHOR{Manasa Mandava}
\AFF{Indian School of Business, Hyderabad,  500032, India}
\AUTHOR{Albert N. Shiryaev}
\AFF{Steklov Mathematical Institute, Moscow,
  119991, Russia, \EMAIL{albertsh@mi-ras.ru}}
} 


\ABSTRACT{This paper extends to Continuous-Time Jump Markov Decision Processes (CTJMDP) the classic result for Markov Decision Processes  stating that, for a given initial state distribution, for every policy there is a (randomized) Markov policy, which can be defined in a natural way, such that at each time instance the marginal distributions of state-action pairs for these two policies coincide.  It is shown in this paper that this equality takes place for a CTJMDP if the corresponding Markov policy defines a nonexplosive jump Markov process.  If this Markov process is explosive,  then at each time instance the marginal probability, that a
state-action pair belongs to a measurable set of  state-action pairs, is not greater for the described Markov
policy  than the same probability for the original policy.  These results are used in this paper to prove that for expected discounted total costs and for average costs per unit time, for a given initial state distribution, for each policy for a CTJMDP the described a Markov policy has the same or better performance.
}

\MSCCLASS{Primary: 90C40; secondary: 90C39, 60J25}

\maketitle
%
\section{Introduction.}\label{S-I}
One of the fundamental facts in the theory of discrete-time Markov decision processes (MDPs) states that, if an initial state distribution is fixed, then for an arbitrary policy there exists a Markov policy such that for each time epoch the marginal distributions of state-action pairs for these policies coincide. This fact was established by \citet{DS} and \citet{Str}. The corresponding Markov policy chooses actions  at a state $x$ at an epoch $t$ with the same probability distribution as the original policy does at epoch $t$ under the condition that the process visits the state $x$ at epoch $t.$  Since most of the major performance criteria, including  expected total discounted costs and average costs per unit time, depend only on marginal distributions of state-action pairs, this theorem implies that for a given initial state distribution  the optimal values of major performance criteria for the classes of all and Markov policies are equal. This fact means that, if an initial distribution is fixed, the decision maker may use only (possibly randomized) Markov policies.

This paper establishes the similar fact for continuous-time jump Markov decision processes (CTJMDPs).  We show that the described equality of state-action marginal distributions holds if the jump Markov process defined by the corresponding Markov policy is nonexplosive, that is, the number of jumps on each finite interval of time is finite with probability 1.  Of course, the notion of explosiveness is applicable only to continuous-time processes because stochastic sequences are always nonexplosive.  In general, when the corresponding Markov policy may define an explosive  Markov process, the results of this paper imply  an inequality rather than the equality. According to this inequality,  at each time instance the probability of the event, that a state-action pair belongs to a measurable set of  state-action pairs for the described Markov policy, is not greater than the probability of the same event for the original policy.   The question, whether the equality always holds, remains open.

We recall that for discrete time the equivalence of arbitrary and Markov policies is proved by induction in the time parameter in \citep{DS, Str} and in other places.  Of course, induction in the time parameter is not applicable to continuous time problems.  Our proofs are based on the fact that for CTJMDPs transition probabilities for marginal distributions of states satisfy an analogue to Kolmogorov's forward equation. This fact was established by \citet{Kit} for problems with bounded intensities, strengthened by \citet{GS} and by \citet{PZ1} to certain nonexplosive CTJMDPs, and further strengthened in this paper to possibly explosive CTJMDPs. 

Studies of Kolmogorov's equations for jump Markov processes were pioneered by \citet{Fel},  who investigated  processes with possibly unbounded transition rates continuous in the time parameter.  In  particular, reference \citep{Fel} includes  the results stating that transition probabilities of  jump Markov processes are solutions of   Kolmogorov's forward equation.  Feller clarified  later in the erratum to \citep{Fel} that these results were obtained  in \citep{Fel}  only for nonexplosive processes. The question whether transition probabilities for explosive jump Markov processes are minimal solutions of Kolmogorov's forward equation remained open for a long time. This question  and the similar question for Kolmogorov's backward equation were answered positively in  Feinberg et al.  \cite{FMS}, where  jump Markov process with
measurable transition rates were considered; see also additional results in \citep{FMS1}.  This fact for forward equation plays the central role in the proofs in this paper.

Section~\ref{S-MD} of this paper describes the model.  Section~\ref{S-MR} presents the main results whose proofs are provided in Section~\ref{S-PM}.  Section~\ref{S-FE} applies the results on  Kolmogorov's forward  equations from \cite{FMS,FMS1} to Markov processes defined by Markov policies for CTJMDPs.   In Section~\ref{S-ER} we show that Markov policies are not worse than arbitrary ones for the total expected discounted  costs  and for average-cost per unit time.

We remark that for discrete-time MDPs, a general policy makes decisions based on the history that includes the current state, time, and previous states and actions.  The decisions may be randomized.  Defining past-dependent and randomized policies for CTJMDPs is a more delicate matter than for discrete-time MDPs.  Early publications on CTJMDPs \cite{GHL,Kak,Mil1,Mil2} dealt only with Markov policies.  For Markov policies the corresponding Markov processes are defined via  Kolmogorov's forward equation; see Guo and Hern\'andez-Lerma~\cite{GHL} for details, where the theory of countable-state CTJMDPs is described. In particular, the results of Section~\ref{S-ER} imply that many currently available results on the existence of optimal and nearly optimal policies within the class of randomized Markov policies hold in the stronger sense because the corresponding policies are also optimal or nearly optimal within the broader class of history-dependent policies. 

\citet{Yus,Yus1} defined history dependent policies for CTJMDPs with bounded jump rates.  For such policies, decisions can be chosen at any time and they depend on the finite sequence $x_0,t_1,x_1,\ldots,t_{n},x_n,t$  of the previous states and jump epochs, the current state, and the current time, where  $0<t_1<t_2<t_3<,\ldots$ and  $t_n<t.$  The arguments relevant to the  Ionescu Tulcea theorem~\cite[Proposition V.1.1]{Neveu} were used in \cite{Yus,Yus1} to construct stochastic processes defined by policies and initial state distributions. These arguments can be extended to CTJMDPs with unbounded jump rates.
\citet{Yus2} also used this approach to reduce CTJMDPs with expected total costs to discrete-time MDPs with actions being the functions of the time parameters with values in the action space.  \citet{Fei94} described the relations between the expected times, during which actions are used between jumps, and actions chosen at the jump epochs.  \citet{Fei04, Fei12} used this relation to reduce discounted CTJMDPs to MDPs with the same action sets as in the CTJMDP.

\citet{Kit} used Jacod's \cite{Jac} results on multivariate point processes and their compensators to provide an equivalent and more elegant construction of stochastic processes for possibly history dependent policies for CTJMDPs.  It was observed by \citet{Kit} that, for a properly defined sample space $(\Omega,\mathcal{F})$ with a filtration, each policy explicitly defines a predictable random measure.  This random measure and an initial state distribution define a jump stochastic process such that the compensator of the random measure for the multivariate point process generated by jumps of this stochastic process is the predictable random measure defined by the policy.  Thus, a policy and an initial state distribution define the appropriate jump stochastic process.

In conclusion, we remark that the notions of randomized policies have different meanings for MDPs and CTJMDPs.  For MDPs, a randomized policy may choose actions randomly at each time instance. For continuous time such policies may not define measurable stochastic processes; \citet[Example 1.2.5]{Kal}. For CTJMDPs, randomized policies are defined as regular policies for the problem with action sets replaced with the sets of probability measures on action sets.  For example, this means that randomized policies may use  transition rates being convex combinations of  transition rates in the original models.  For this reason, randomized policies for CTJMDPs are often called relaxed, which is a more precise term. In this paper we mostly consider only relaxed policies, and the terms ``randomized" and ``relaxed" are  used only at the end of Section~\ref{S-ER}.

\section{Model description.}\label{S-MD}
 In this section we introduce basic notations, define CTJMDPs,  and provide the construction of jump stochastic processes defined by policies and initial state distributions. Recall that a measurable space $(S, \mathcal{S})$  is called a standard Borel space, if there is a measurable one-to-one correspondence $f$ of this space onto a Polish (complete, separable, metric) space endowed with its Borel $\sigma$-algebra such that the correspondence $f^{-1}$ is also measurable. We usually write $(S, \B(S))$ instead of $(S, \mathcal{S})$ for a standard Borel space. If $S'\in \B(S),$ then we consider the $\sigma$-algebra $\B(S')=\{\tilde{S}\cap S^\prime:\,\tilde{S} \in\B(S)\}$ on the set $S'.$  The measurable space $(S',\B(S'))$ is also a standard Borel space. If we add an isolated point $s' \notin S$ to a standard Borel space $S$, then we consider the $\sigma$-algebra $\B(S \cup \{s'\}) = \sigma(\{\B(S), \{s'\}\})$, where $\sigma(\mathcal{E})$ denotes the $\sigma$-algebra  generated on a set by the set $\mathcal{E}$ of its subsets. Of course, $(S \cup\{s'\}, \B(S \cup \{s'\})$ is also a standard Borel space. We denote by $\P(S)$  the set of all probability measures on $(S, \B(S)).$ For two standard Borel spaces $(S, \B(S))$ and $(\tilde{S}, \B(\tilde{S})),$ a transition probability $\pi(\cdot|\cdot)$ from $(S, \B(S))$ to $(\tilde{S}, \B(\tilde{S}))$ is a mapping from $(S, \B(S))$ to $\P(\tilde S)$ such that for each $E\in\B(\tilde{S})$ the function $\pi(E|s) :(S,\B(S))\mapsto ([0,1],\B([0,1]))$ is measurable. A Dirac measure concentrated at a point $s$ is denoted by $\delta_s.$  Let  $\R : = ]-\infty, +\infty[,$ $\bar{\R}= [-\infty, +\infty]$, $\R_+ : = ]0, +\infty[$, $\bar{\R}_+ := ]0, +\infty]$, and $\R_+^0:= [0,+\infty[.$

The probability structure of a CTJMDP  is specified by the four objects $\{ X, A, A(\cdot), \tilde{q} \}$, where
\begin{itemize}
\item[(i)] $(X, \B(X))$ is a standard Borel space (the state space);
\item[(ii)] $(A, \B(A))$ is a standard Borel space (the action space);
\item[(iii)] $A(x)$ is a non-empty subset of $A$ for each state $x \in X$ (the set of actions available at $x$). It is assumed that the set of feasible state-action pairs \[Gr(A): = \{ (x,a) : x \in X, a \in A(x)\}\] is a measurable subset of $(X \times A)$ containing the graph of a measurable mapping of $ X$ to $A$.
\item[(iv)] $\tilde{q}( x, a, \cdot)$ is  a signed measure on $(X,\B(X))$ for each $(x,a)  \in Gr(A)$ (the transition rate), such that $\tilde{q}(x,a, X) = 0$, $0 \le \tilde{q}(x,a, Z \setminus \{x\})  < \infty$, and $\tilde{q}(x,a,Z)$ is a measurable function on $Gr(A)$ for each $Z \in \B(X)$.
\end{itemize}
Let $\tilde{q}(x,a) := \tilde{q}( x,a,X \setminus \{x\} )$ for all $(x,a) \in Gr(A)$ and let ${\bar q}(x):= \sup_{a \in A(x)}\tilde{q}(x,a)$ for all $x \in X$. If an action $a\in A(x)$ is selected at state $x\in X$ and is fixed until the next jump, then the sojourn time has an exponential distribution with the intensity $\tilde{q}(x,a)$ and the process jumps to the set $Z\setminus\{x\},$  where $Z\in\B(X),$  with probability $\tilde{q}(x,a,Z\setminus\{x\})/\tilde{q}(x,a)$ if $\tilde{q}(x,a)>0.$  If $\tilde{q}(x,a)=0,$ then the state $x$ is absorbing. However, the model allows changing actions between jumps.  In this paper we make the following standard assumption, which implies that there are no instantaneous jumps.
\begin{Assumption}
\label{A1}
$\bar{q}(x) < \infty$ for each $x \in X$.
\end{Assumption}

To define a sample space, which includes trajectories that have a finite number of jumps over $\BB{R}_+$ and that have an infinite number of jumps over a finite interval of time, we add an additional point $x_\infty \notin X$ to $X$. Let ${\bar X} := X \cup\{x_\infty\}$. 
For a set $\mathcal{H}$ of real-valued functions defined on a common set, let $\sigma(\mathcal{H})=\sigma(\{f^{-1}(\B(\R)):\,f\in\mathcal{H}\})$ denote the sigma field generated by all functions from  $\mathcal{H}.$ 

Let $X\times (\bar{\BB{R}}_+ \times {\bar X})^\infty$ be the set of all sequences $(x_0,t_1, x_1,t_2, x_2, \ldots)$ with $x_0\in X,$ $t_{n}\in\bar{\R}_+,$ and $x_n\in \bar{X}$ for $n =1,2, \ldots\ .$  This set is endowed with the $\sigma$-algebra $\B(X\times (\bar{\BB{R}}_+ \times {\bar X})^\infty)$  defined by the products of the Borel $\sigma$-algebras  $\B{(X)},$  $\B(\bar{\R}_+),$ and $\B(\bar{X})$. Since a countable product of standard Borel spaces is a standard Borel space, the measurable space  $(X\times (\bar{\BB{R}}_+ \times {\bar X})^\infty, \B(X\times (\bar{\BB{R}}_+ \times {\bar X})^\infty ))$ is a standard Borel space.

The set $\Omega$ of trajectories with a finite or countable numbers of jumps is defined as  the subset of all sequences $(x_0,t_1,
x_1, t_2, x_2, \ldots)$ from $X\times (\bar
{\BB{R}}_+  \times {\bar X})^\infty  $ such that, for  $n = 1,2,\ldots,$ the  following two properties hold: (i) if $t_n < \infty,$ then $x_n \in X$ and $t_{n+1}> t_n,$   and (ii) if $t_n = +\infty,$ then $x_n = x_\infty$ and $t_{n+1} = +\infty.$  The definition of $\Omega$ implies that $\Omega\in \B(X\times (\bar{\BB{R}}_+ \times {\bar X})^\infty).$  Let us denote $\F=\B(\Omega)=\{\Omega\cap B:\, B\in \B(X\times (\bar{\BB{R}}_+ \times {\bar X})^\infty)\},$ where the second equality is the definition of $\B(\Omega).$  The standard Borel space $(\Omega,\F)$ is called the sample space.

Define the random variables    $t_0(\omega) := 0$, $x_0(\omega) := x_0$, $t_n(\omega):=t_n,$ and  $x_{n}(\omega):=x_{n},$ for  $n = 1,2,\ldots,$ on the  measurable space $(\Omega, \mathscr{F})$ denoting, respectively, the initial  time epoch, the initial state,  the time of the $n$th jump, and the state to which the process jumps at the $n$th jump. Let $t_\infty(\omega) := \lim_{n \to \infty} t_n (\omega).$ 
The jump process of interest, $\{\xi_t(\omega): t \in \BB{R}_+^0, \omega\in\Omega\}$ with values in $\bar{X}$, 
is 
\begin{equation}
\label{JMP}
\xi_t(\omega) = x_n(\omega), \quad \text{ for } t_n(\omega) \leq t < t_{n+1}(\omega), n=0,1,\ldots, \quad \text{and} \quad \xi_t(\omega) =x_\infty  \text{ for } t\ge t_\infty(\omega).
\end{equation}
Let us consider the natural filtration $\F_t=\sigma(\{\xi_s(\omega): 0\le s\le t \}),$ $t\in \BB{R}_+^0,$ and  $\sigma$-algebras $\F_\infty=\sigma(\{{\F_t: t\ge 0}\}) $ and $\F_{t_n}=\{B\in\F:B\cap\{t_n\le t\}\in\F_t, t\ge 0\}.$ The definition of $(\Omega,\F)$ implies that $\F_{t_n}=\sigma(\{t_m,x_m:\, 0\le m\le n\});$ see e.g., \cite[Theorem 4.13]{KR}.


  A policy $\pi$ is a mapping $(\Omega\times\R_+)\mapsto  \P(A)$ such that (i) the stochastic process
  $\pi(B|\omega,t)$ is predictable for all $B\in\B(A),$ and (ii) $\pi(A( \xi_{t^-}(\omega))|\omega,t)=1$ for all $(\omega,t)\in  (\Omega\times\BB{R}_+)$ with  $t< t_\infty(\omega)$.  Since  it is possible that $t_\infty(\omega)<+\infty$   for some $\omega,$ in order to define a policy for all $t\in \BB{R}_+,$ including $t\ge t_\infty(\omega),$ we add an additional point $a_\infty \notin A$ to $A$ and set $A(x_\infty) =  \{a_\infty\}$. Let $\bar{A} := A \cup \{a_\infty\}$. The definition and the structure of predictable processes described in \citet[p. 241]{Jac} implies that $\pi$ is a policy if and only if there is a sequence of
transition probabilities $\pi^n:((X\times \BB{R}_+)^{n+1}, \B((X\times \R_+)^{n+1}))\mapsto(A, \B(A))$ such that, at each  $t\in\BB{R}_+,$ the policy $\pi$ selects the probability measure
 \begin{equation}
\label{policy}
\pi(\,\cdot\,| \omega, t) := \sum_{n \ge 0}\pi^n(\,\cdot\,| x_0, t_1, x_1, \ldots, t_n, x_n, t-t_n)I\{ t_n < t \le t_{n+1}\} + \delta_{a_\infty} (\cdot) I\{t \ge t_\infty\}, \ \  \omega \in \Omega,  
\end{equation}
where we omit $\omega$ in the right-hand side of \eqref{policy} and $\delta_{a_\infty}(\cdot)$ is a Dirac measure on $(\bar{A}, \B(\bar{A}))$ concentrated at $a_\infty.$ 

 A policy $\pi$ is called {\it Markov} if there exists  a  transition probability $\tilde{\pi}$ from $((X\times \BB{R}_+), \B((X\times \R_+)))$ to $(A, \B(A))$ such that
 $\pi(\cdot|\omega,t)=\tilde{\pi}(\cdot|\xi_{t-}(\omega),t)$ for all $(\omega,t)\in  (\Omega\times\BB{R}_+)$ with  $t< t_\infty(\omega)$. For a Markov policy $\pi$, formula \eqref{policy} implies that $\pi^n(B|x_0, t_1, x_1, \ldots, t_n, x_n, t-t_n)=\tilde{\pi}(B|x_n,t),$   when $t_n < t \le t_{n+1}$ and   for all $B\in\B(X)$ and  $n= 0,1,2,\ldots\ .$  With a slight abuse of notations, we shall write $\pi$ instead of $\tilde{\pi}.$

For a measurable function $f$ on $X \times A$, define
\begin{equation}
\label{Ext}
f(z,p) := \int_{A(z)} f(z,a)p(da), \qquad z \in X,\ p \in \P(A),
\end{equation}
whenever the integral is defined. In particular, \eqref{Ext} for $f(z,a) = \tilde{q}(z,a,Z)$ is
\begin{equation}
\label{q-def21}
\tilde{q}( z, p, Z) =
\int_{A(z)} \tilde{q}(z, a, Z) p (da), \qquad z \in X,\ p \in \P(A),\ Z \in \B(X),
\end{equation}
and, \eqref{Ext} for $f(z,a) = \tilde{q}(z,a)$ is
\begin{equation}
\label{q-def22}
\tilde{q}(z,p) = \int_{A(z)}\tilde{q}(z,a) p(da), \qquad z \in X,\ p \in \P(A),
\end{equation}
where we set $\tilde{q}(z,a):=0$ for $(z,a)\in (X\times A)\setminus Gr(A);$ in particular $ \tilde{q}(z, a, Z)=0$ if $(z,a)\in (X\times A)\setminus Gr(A)$ and $Z\in\B(X).$
Due to Assumption~\ref{A1}, the integrals in \eqref{q-def21} and \eqref{q-def22} are defined, and
\begin{equation}
\label{f1}
\tilde{q}(z, p )  \le \int_{A(z)} \left(\sup_{a \in A(z)} \tilde{q}(z,a)\right) p(da) \le  \bar{q}(z) < +\infty, \qquad z \in X,\  p \in \P(A).
\end{equation}
In addition, the properties of the transition rate $\tilde{q}(z,a,Z)$  imply that $\tilde{q}(z,p, Z)$ is a signed measure on $(X,\B(X))$ with $\tilde{q}(z,p,X) = 0$, the function $\tilde{q}(z,p,Z\setminus\{z\})$ is a finite measure on $(X,\B(X))$. For a policy $\pi$, let  $\pi_t(\omega)$  denote the probability measure with values  $\pi(\,\cdot\, | \omega, t).$  Then, $q(\xi_{t-}(\omega), \pi_t(\omega), Z\setminus \{\xi_{t-}(\omega)\})$ defined by \eqref{q-def21}, with $z = \xi_{t-}(\omega)$, $p(\cdot) = \pi_{t}(\omega)$, and $Z = Z \setminus\{\xi_{t-}(\omega)\}$, is the jump intensity at time $t$ from the state $\xi_{t-}(\omega)$ to the  set $Z\setminus\{\xi_{t-}(\omega)\}.$


Recall that a multivariate point process is a sequence $(t_n(\omega), x_n(\omega))_{n\ge 1}$ of random variables on $(\Omega,\F)$ with values in $(\BB{\bar{R}}_+ \times \bar{X})$ and such that, for  $n = 1,2,\ldots,$ the following properties hold: (i) if $t_n < +\infty,$ then $x_n \in X$ and $t_{n+1}> t_n,$   and  if $t_n = +\infty,$ then $x_n = x_\infty$ and $t_{n+1} = +\infty,$ (ii) $t_n(\omega)$ is a stopping time, and (iii) $x_n(\omega)$ is $\F_{t_n}-$measurable. A multivariate point process $(t_n(\omega), x_n(\omega))_{n\ge 1}$ is characterized by the  random measure $\mu$ on $({R}_+^0 \times X)$ defined by
\begin{equation} \label{m1}
\mu(\omega; [ 0,t], Z) := \sum_{ n \geq 1} I\{t_n(\omega) \in \,[ 0,t]\} I\{x_n(\omega) \in Z\}, \qquad  \omega \in \Omega,\  t \in \mathbb{R}_+^0  ,\ Z \in \B(X).
\end{equation}

A random measure $\nu$ on $\R_+\times X$ is called predictable if for every $Z\in\B(X)$ the stochastic process $\{\nu(\omega; [0,t], Z)\}$ is $\F_{t-}$-measurable.  According to \citet[Theorem~2.1]{Jac} or \citet[Theorem~4.20]{KR}, for a given probability space $(\Omega, \F, \BB{P})$ with a right-continuous filtration $\{\F_t\}_{t \ge 0},$ there exists a predictable random measure $\nu:(\Omega \times \B(\BB{R}_+^0 \times X)) \to \BB{R}_+^0$ called the {\it compensator} of $\mu$ such that (i) for each $Z \in \B(X)$, the  process $\{\nu(\omega; [0,t], Z)\}_{t  \in \BB{R}_+^0}$ is predictable; and (ii) for any stopping time $T$ with values in $\BB{R}_+$ and $Z \in \B(X),$
 \begin{equation}
\label{mu-nu}
\BB{E}(\mu(\omega;[0,T], Z)) = \BB{E} (\nu(\omega;[0,T], Z)),
\end{equation}
where $\BB{E}$ denotes the expectation with respect to the probability measure $\BB{P}.$

%
Define the random measure $\nu^\pi$ on $({R}_+^0 \times X)$ by
\begin{equation}
\label{m2}
\nu^\pi (\omega; [ 0,t], Z) := \int_{0}^{t} \tilde{q}( \xi_{s}(\omega), \pi_s(\omega), Z \setminus \{\xi_{s}(\omega)\})I\{\xi_{s}(\omega) \in X\} ds, \quad  \omega \in \Omega,\  t \in \mathbb{R}_+^0,\ Z \in \B(X).
\end{equation}
This random measure is predictable.  Indeed, in view of \eqref{JMP} and \eqref{policy}, for each $Z\in\B(X),$ the stochastic process $\{\nu^\pi(\omega; [0,t], Z)\}$ is $\F_{t}$-measurable. In addition, it has continuous paths.  Therefore, these processes are $\F_{t-}$ measurable or, in other words, predictable; see, e.g., \citet[Proposition~2.6]{JS} or \citet[Theorem 4.16]{KR}.

Furthermore, $\nu^\pi(\omega; [t_\infty, +\infty[, X) = 0$ since $\xi_t(\omega) = x_\infty$ for all $t \ge t_\infty$ and $\nu^\pi(\omega; \{t\} \times X) = 0$ since the function $\nu^\pi(\omega; [ 0,t], X)$ is continuous in $t\in \R_+.$ In view of \citet[Theorem~3.6]{Jac}, the predictable random measure $\nu^\pi$   and a probability measure $\gamma$ on $X$ define a unique probability measure $\BB{P}_\gamma^\pi$ on $(\Omega,\F)$ for which $\BB{P}_\gamma^\pi(dx_0)=\gamma(dx_0)$ and  $\nu^\pi$ is a compensator of the random measure   $\mu.$  We remark that \cite[Theorem~3.6]{Jac} has two assumptions, namely, \cite[assumptions (4) and (A.2)]{Jac}.  Assumption (4) from \cite{Jac} is verified in the first sentence of this paragraph. Assumption (A.2) follows from the construction of the sample space $(\Omega,\F).$

 If $\gamma(\{x\}) = 1$ for some $x \in X$, we shall write $\BB{P}_x^\pi$ instead of $\BB{P}_\gamma^\pi$.   Let $\BB{E}_\gamma^\pi$ and $\BB{E}_x^\pi$  denote  expectations with respect to the measures  $\BB{P}_\gamma^\pi$ and $\BB{P}_x^\pi$ accordingly.  For a policy $\pi$ and an initial distribution $\gamma$, we say that the jump process is \textit{nonexplosive} if $\BB{P}_\gamma^\pi(\xi_t(\omega) \in X) = 1$ for all $t \in \BB{R}_+$.

\citet{Yus,Yus1} constructed explicitly the probability measure $\BB{P}_\gamma^\pi$ for a given initial distribution $\gamma$ and a nonrandomized policy $\pi$ by using the Ionescu Tulcea theorem.   The policies we consider in this paper can be viewed as nonrandomized if  actions are substituted with probability measures on the feasible sets of actions.  Therefore, in view of the Ionescu Tulcea theorem the function $\BB{P}_x^\pi(C)$ is measurable in $x$ for every $C\in \F.$


 Observe that, for all $t \in \BB{R}_+$,
\begin{equation}
\label{zero-m}\BB{P}_\gamma^\pi(\xi_{t}(\omega) \not = \xi_{t-}(\omega)) = \BB{E}_\gamma^\pi \mu(\omega; \{t\}, X)  = \BB{E}_\gamma^\pi \nu^\pi(\omega;\{t\},X) = 0,\end{equation}
where the first equality follows from the definition of a random measure of a multivariate point process, the second equality follows from \eqref{mu-nu}, and the last one is follows from $\nu^\pi(\omega;\{t\},X)=0 .$ 

We now define in \eqref{MM} marginal distributions of the states and of the state-action pairs. Consider the process of actions $\{U_t(\omega): t \in \BB{R}_+, \omega \in \Omega\}$ with values in $\bar{A},$ where the probability of $U_t(\omega) \in B$ is $\pi(B|\omega, t)$ for $t<t_\infty(\omega)$ and $U_t(\omega) = a_\infty$ for $t \ge t_\infty(\omega)$.
For the given $\omega\in\Omega$ and $t>0,$ the probability of the event $\{\xi_{t-}(\omega)\in Z, U_t(\omega)\in B \},$  where $ Z \in \B(\bar{X})$ and  $B \in \B(\bar{A}),$   is $I \{ \xi_{t-}(\omega) \in Z\}\pi(B|\omega,t).$

 For an initial distribution $\gamma$ and a policy $\pi$, consider the marginal probabilities
\begin{align}
\label{mm}
P_\gamma^{\pi} (t, Z) &:= \BB{P}_\gamma^\pi ( \xi_t(\omega) \in Z)= \BB{P}_\gamma^\pi ( \xi_{t-}(\omega) \in Z),\\
\label{MM}
P_\gamma^{\pi} (t, Z, B) &:= \BB{P}_\gamma^\pi ( \xi_{t-}(\omega) \in Z,U_t(\omega)\in B)= \BB{P}_\gamma^\pi ( \xi_t(\omega) \in Z,U_t(\omega)\in B).
\end{align}
where $t \in \BB{R}_+$, $Z \in \B(\bar{X})$, and $B \in \B(\bar{A})$. The second equalities in \eqref{mm} and \eqref{MM} are correct in view of \eqref{zero-m}.


Observe that, in view of \eqref{zero-m}--\eqref{MM}, for $ t \in \BB{R}_+,$ $ Z \in \B(\bar{X}),$ and  $B \in \B(\bar{A}),$
\begin{eqnarray}
\label{m}
P_\gamma^\pi(t, Z) &=&  P_\gamma^\pi(t,Z,\bar{A}), \\
\label{M}
P_\gamma^\pi(t, Z,B) &=&  \BB{E}_\gamma^\pi (I \{ \xi_{t-}(\omega) \in Z\}\pi(B|\omega,t))= \BB{E}_\gamma^\pi [I \{ \xi_{t}(\omega) \in Z\}\pi(B|\omega,t)],
\end{eqnarray}
and the function $P_\gamma^\pi(\cdot, \cdot, \cdot)$ is a transition probability from $(\BB{R}_+, \B(\BB{R}_+))$ to $(\bar{X}\times\bar{A},\B(\bar{X}\times\bar{A}))$. 
To see that $P_\gamma^\pi(\cdot, \cdot, \cdot)$ is a transition probability from $(\BB{R}_+, \B(\BB{R}_+))$ to $(\bar{X}\times\bar{A},\B(\bar{X}\times\bar{A}))$, observe that by its definition $P_\gamma^\pi(t, \cdot, \cdot)$ is a probability measure on $(\bar X \times \bar A)$ for all $t\in\R_+$.  Since the processes $\xi_t(\omega)$ and $\pi(U_t(\omega)\in B |  \omega,t)$  are measurable and predictable, respectively,  the processes $I\{\xi_t(\omega)\in Z\}$ and $\pi(U_t(\omega)\in B |  \omega,t)$  are measurable.  
Since \eqref{M} holds and the functions under the expectation in \eqref{M} are measurable in $(\omega,t),$  the function $P_\gamma^\pi(t,Z,B)$ is measurable in $t.$ In addition, for any nonnegative measurable function $f$ on $Gr(A)$,
\begin{multline}
\label{f-Ex}
\resizebox{.9\hsize}{!}{$\BB{E}_\gamma^\pi f(\xi_{t}(\omega), \pi_t(\omega))I\{\xi_{t}(\omega) \in Z\} = \BB{E}_\gamma^\pi I\{\xi_{t}(\omega) \in Z\}\int_{A} \Big( f(\xi_{t}(\omega), a) I\{a \in A(\xi_{t}(\omega))\} \Big) \pi(da \mid \omega, t) $}\\
 = \int_Z \int_{A} f(z,a)\delta_{a}(A(z)) P_\gamma^\pi(t,dz,da) = \int_Z \int_{A(z)} f(z,a)P_\gamma^\pi(t,dz,da),
\end{multline}
where the first equality follows from \eqref{Ext}, the second equality follows from \eqref{M}, and the last one is straightforward. Similar to $\BB{P}_x^\pi$, we shall write $P_x^{\pi}$ instead of $P_\gamma^{\pi}$  if $\gamma(\{x\}) = 1$ for some $x \in X$. In the rest of this paper, we omit $\omega$ whenever there is no confusion.

\section{Main results.}
\label{S-MR}
In this section we formulate the main result of this paper. Let us fix an arbitrary $B\in\B(\bar{A})$ and consider the measures $P^\pi_\gamma(t,\cdot,B)$ and $P^\pi_\gamma(t,\cdot)$ on $(\bar{X},\B(\bar{X})),$ where $t\in\mathbb{R}_+.$ Then, \eqref{m} implies $P^\pi_\gamma(t,Z,B)\le P^\pi_\gamma(t,Z)$ for all $Z\in \B(\bar{X}).$ Thus $P^\pi_\gamma(t,\cdot,B)\ll P^\pi_\gamma(t,\cdot).$ Therefore, in view of the Radon-Nikodym theorem, there is a derivative $\frac{dP^\pi_\gamma(t,\cdot,B)}{dP^\pi_\gamma(t,\cdot)}.$  The following lemma and its corollary state that this derivative can be written as a Markov policy.
\begin{lemma}
\label{def}
For an initial distribution $\gamma$ on $X$ and a policy $\pi$, there exists a Markov policy $\varphi$ such that, for all $t \in \BB{R}_+,$ 
\begin{equation}
\label{formula-de2018}
P_\gamma^{\pi} (t, Z, B) =\int_Z \varphi(B \lvert z, t)P_\gamma^{\pi}(t, dz), \quad Z\in \B(X),\ B \in \B(A).
%
\end{equation}
%
\end{lemma}
\proof{Proof.}
As explained after formula~\eqref{M}, the function $P_\gamma^\pi(\cdot, \cdot, \cdot)$ is a transition probability from $(\BB{R}_+, \B(\BB{R}_+))$  to $(\bar{X}\times\bar{A},\B(\bar{X}\times\bar{A}))$. Therefore, in view of \eqref{m}, there exists a transition probability $\tilde{\varphi}$ from $(\bar{X} \times \BB{R}_+,\B(\bar{X} \times \BB{R}_+))$ to $(\bar{A},\B(\bar{A}))$ such that formula~\eqref{formula-de2018} holds with $\varphi = \tilde{\varphi}$; see e.g., Bertsekas and Shreve~\cite[Corollary 7.27.1]{BS}. In addition, since the action chosen by the policy $\pi$ at time $t$ is concentrated on $A(\xi_{t-})$,  the probability measure $P_\gamma^\pi(t, \cdot, \cdot)$ is concentrated on $\{Gr(A) \cup (x_\infty, a_\infty)\}.$ Thus, for all $t \in \BB{R}_+$, 
\begin{equation}
\label{a-s}
\tilde{\varphi}(A(z)|z,t) = 1, \quad z \in X\  (P_\gamma^\pi(t, \cdot)\emph{-a.s.}).
\end{equation}

Suppose that $\phi$ is a measurable mapping from $\bar{X}$ to $\bar{A}$ with $\phi(x) \in A(x)$ for all $x \in \bar{X}$. The existence of such a mapping is guaranteed by assumption (iii) in Section~\ref{S-MD} stating that the set of feasible state action pairs $Gr(A)$ contains the graph of a measurable mapping $\phi$ from $X$ to $A$.
Then, the function $\varphi(B|z,t)$, where $B \in \B(A)$, $z \in X$, and $t \in \BB{R}_+$, defined by
\begin{equation}
\label{varphi}
\varphi(B|z,t) := \left\{ \begin{array}{ll}
\tilde{\varphi}(B|z,t), \quad &\text{ if } \tilde{\varphi}(A(z)|z,t) =1,\\
\delta_{\phi(z)}(B), \quad &\text{ otherwise, }
\end{array}
\right.
\end{equation}
is a measure on $(A,\B(A))$ with $\varphi(A(z)|z,t) = 1$ for each $z \in X$ and $t \in \BB{R}_+$. In addition, for every $B \in \B(A)$, the function $\varphi(B\mid z,t)$ is measurable in $(z,t)$. To see this, observe that
\[\tilde{\varphi}(A(z)|z,t)=  \int_A \delta_{a}(A(z))\tilde{\varphi}(da |z,t), \quad z \in X,\ t \in \BB{R}_+.\]
Since the set $Gr(A) \in \B(X \times A)$, the function $\delta_{a}(A(z))$ is measurable on $X \times A.$  The measurable property of integrals with respect to a transition probability imply that the function $\tilde{\varphi}(A(z)|z,t)$ is measurable in $(z,t)$; see e.g., \cite[Proposition~2.9]{BS}. Then, the set $\{(z,t):\tilde{\varphi}(A(z)\mid z,t) = 1 \}$ is measurable. This fact implies that  the function $\varphi(B\mid z,t)$ is measurable in $(z,t)$ for every $B \in \B(A)$. Hence, the function $\varphi$ defined by \eqref{varphi} is a Markov policy. Therefore, in view of \eqref{a-s}, it follows from \eqref{formula-de2018} with $\varphi=\tilde{\varphi}$ that \eqref{formula-de2018} holds for the Markov policy $\varphi$ defined by \eqref{varphi}.
\hfill \Halmos\endproof

\begin{remark}
Strictly speaking, a Markov policy $\varphi$ satisfying \eqref{formula-de2018} depends on the initial distribution $\gamma$ and on the policy $\pi$, that is, $\varphi = \varphi_{\gamma,\pi}$. To simplify notations, we write $\varphi$ instead of $\varphi_{\gamma,\pi}$. 
\end{remark}

\begin{corollary}\label{cordef}
For an initial state distribution $\gamma$  on $X$ and for a policy $\pi,$ a Markov policy $\varphi,$ whose  existence is stated in Lemma~\ref{def}, satisfies \eqref{formula-de2018} if and only if, for all   $t \in \BB{R}_+$  and $B \in \B(A),$ 
\begin{equation}
\label{formula-d}
\varphi(B \lvert z, t) = \frac{P_\gamma^{\pi}(t, dz,B)}{P_\gamma^{\pi}(t, dz)}, \quad   z \in X\  (P_\gamma^\pi(t, \cdot)\it{-a.s.}).
\end{equation}
\end{corollary}

\proof{Proof.} The  corollary follows from the definition of the Radon-Nikodym derivative.
\hfill \Halmos\endproof

The following theorem is the main result of this paper.
\begin{theorem}
\label{main}
For an initial distribution $\gamma$ on $X$ and a policy $\pi$, let $\varphi$ be a Markov policy satisfying \eqref{formula-de2018}. Then
\begin{equation}
\label{2-d}
P_\gamma^{\varphi} (t, Z, B) \le P_\gamma^{\pi} (t, Z, B), \qquad t \in \mathbb{R}_+, Z \in \B(X), B \in \B(A).
\end{equation}
In addition, if $P_\gamma^\varphi(s,X) = 1$ for some $s \in \BB{R}_+$, then \eqref{2-d} holds for all $t\in ]0,s]$ with an equality.  In particular, if $P_\gamma^\varphi(t,X) = 1$ for all $t \in \BB{R}_+$, then \eqref{2-d} holds  with an equality.
\end{theorem}
%
\begin{corollary}
\label{C:main}
Let the transition rates $q(z,a)$ be bounded in  $(z,a)\in Gr(A)$. Then, for every policy $\pi$ and initial distribution $\gamma$, $P_\gamma^\pi(t,X) = 1$ for all $t \in \BB{R}_+.$ In addition, formula~\eqref{2-d} holds with an equality for every Markov policy $\varphi$ satisfying \eqref{formula-de2018}.
\end{corollary}
\proof{Proof.}
Let $\lambda$ be a non-negative integer, such that $q(z,a) < \lambda$ for all $(z,a) \in Gr(A)$, and $N(t)$ be a Poisson process with  the rate $\lambda$.  Then, $P_\gamma^\pi(t_\infty < \infty)  \le \BB{P}(N(t) = \infty \ for \ some \ t \in \BB{R}_+) = 0.$  In view of this fact, \eqref{mm}, and  $\{t< t_\infty\} = \{\xi_t \in X\}$,
\begin{equation}
\label{nonexp}
P_\gamma^\pi(t,X) = \BB{P}_\gamma^\pi(\xi_t \in X) = \BB{P}_\gamma^\pi(t < t_\infty) = 1, \qquad t \in \BB{R}_+.
\end{equation}
The second statement in the corollary follows from Theorem~\ref{main} and \eqref{nonexp}  applied to the Markov policy satisfying \eqref{formula-de2018}.
 \hfill\Halmos\endproof
\begin{corollary}
For an initial distribution $\gamma$ on $X$ and a policy $\pi$, let $\varphi_1$ and $\varphi_2$ be two Markov policies satisfying \eqref{formula-de2018}. Then $\BB{P}_\gamma^{\varphi_1}=\BB{P}_\gamma^{\varphi_2}$ and
\begin{equation}
\label{2-deuf}
P_\gamma^{\varphi_1} (t, Z, B) = P_\gamma^{\varphi_2} (t, Z, B), \qquad t \in \mathbb{R}_+,  Z \in \B(X), B \in \B(A).
\end{equation}
\end{corollary}
\proof{Proof.}
Theorem~\ref{main} implies that $P_\gamma^{\varphi_i} (t, \cdot) \ll P_\gamma^{\pi} (t, \cdot)$ for all $t\in \BB{R}_+$ and for $i=1,2.$  Therefore, in view of \eqref{formula-d},
$\varphi_1(B|x,t)=\varphi_2(B|x,t)$ $P_\gamma^{\varphi_1}(t,\cdot)$--a.s. and $P_\gamma^{\varphi_2}(t,\cdot)$--a.s.,   $z\in X,$ for all $t\in \BB{R}_+$ and $B \in \B(A).$
In view of Theorem~\ref{main} applied to  $\varphi=\varphi_1$ and $\pi=\varphi_2,$ the inequality
$
P_\gamma^{\varphi_1} (t, Z, B) \le P_\gamma^{\varphi_2} (t, Z, B)$ holds for all $ t \in \mathbb{R}_+,$ $  Z \in \B(X),$ and $ B \in \B(A).$
The same arguments imply that $P_\gamma^{\varphi_2} (t, Z, B) \le P_\gamma^{\varphi_1} (t, Z, B).$  Thus, \eqref{2-deuf} is proved.  In particular,
$P_\gamma^{\varphi_1} (t, \cdot) = P_\gamma^{\varphi_2} (t,\cdot)$ for all $t\in \BB{R}_+.$  Thus, two jump Markov processes have the same marginal distributions.  This implies that their distributions coincide.
 \hfill\Halmos\endproof

The proof of Theorem~\ref{main} is given in Section~\ref{S-PM}. Theorems~\ref{thm:S1} and \ref{thm:S1inf}  from Section~\ref{S-ER}, which follow from  Theorem~\ref{main}, establish the sufficiency of Markov policies for particular criteria.  In Section~\ref{S-FE} we present  auxiliary facts that follow from the results in Feinberg et al.~\cite{FMS1}.

\section{Kolmogorov's forward equation corresponding to a Markov policy.}
\label{S-FE}

In this section we verify that (i) the jump process corresponding to a Markov policy is a  jump Markov process, and (ii) its transition probability is the minimal solution of  Kolmogorov's forward equation. These facts follow from Feinberg et al.~\cite[Theorem~1 and Corollary~7]{FMS1}. For readers' convenience, we present these two results here. Note that, action sets are not considered in \cite{FMS1}, jump intensities are defined there by $Q$-functions, whose definition is given below, and  jump Markov processes are defined in \cite{FMS1} on a finite or infinite interval $[T_0, T_1[\in \BB{R}_+$. We present here the results for the case  $T_0 = 0$ and $T_1 = \infty$ needed in this paper.

Recall that, a function $q(z,t,Z),$ where $z \in X$, $t \in \BB{R}_+$, and $Z \in  \B(X),$ is called a $Q$-function if:

(a) for fixed $z$  and $t,$ the function $q(z,t,Z)$ is a signed measure on $(X,\B(X))$ with $q(z,t,X) = 0$ and $q(z,t,Z\setminus\{z\})$ is a finite measure on $(X,\B(X));$

(b) for a fixed $Z$ the function $q(z,t,Z)$ is measurable in $(z,t)$.

Let us consider a $Q$-function $q$ satisfying the condition
\begin{equation}
\label{A-4}
\int_0^t q(z,s) ds < \infty, \qquad t \in \BB{R}_+, z \in X.
\end{equation}
Then the $Q$ function $q$ defines the   predictable random measure $\nu$ on $(R_+^0 \times X)$ such that
\begin{equation}
\label{m22}
\nu (\omega; [ 0,t], Z) := \int_{0}^{t} q( \xi_{s}(\omega), s, Z \setminus \{\xi_{s}(\omega)\})I\{\xi_{s}(\omega) \in X\} ds, \qquad  \omega \in \Omega,\  t \in \mathbb{R}_+^0,\ Z \in \B(X),
\end{equation}
where the function $\xi_{s}(\omega)$ is defined in \eqref{JMP}.  An initial state distribution $\gamma$ on $X$ and a compensator $\nu$  of the random measure for the multivariate point process  $(t_n(\omega),x_n(\omega))_{n\ge 1}$  uniquely define a probability distribution $\BB{P}_\gamma$ on $(\Omega,\cal{F})$ such that $\BB{P}_\gamma(x_0\in Z)=\gamma(Z)$ for all $Z\in\cal{B}(X).$ 


Let $q(x,t): = q(x, t, X \setminus\{x\})$ for $x \in X$ and $t \in \BB{R}_+$. Following \citet[Theorem 2]{Fel}, for  $u \in \BB{R}_+^0$, $x \in X$, $t \in ]u,T_1[$, and $Z \in \B( X)$, define
\begin{equation}
\label{b0}
\bar{P}^{{(0)}} (u,x;t,Z) = \delta_{x}(Z) e^{-\int_u^t
q(x, s) ds},
\end{equation}
and, for $n = 1,2,\ldots,$
\begin{equation}
\label{bn}
\bar{P}^{(n)}(u, x; t, Z) = \int_{u}^{t} \int_{ X  }
 e^{ -\int_{u}^{w} q(x,\theta ) d\theta }  q(x,w,dy \setminus \{x\}) \bar{P}^{(n-1)}(w, y; t, Z)
 dw.
\end{equation}
Set
\begin{equation}
\label{edef}
\bar{P}(u, x; t, Z) = \sum\limits_{n=0}^{\infty} \bar{P}^{(n)}(u, x; t, Z).
\end{equation}

\begin{theorem} {\rm (\cite[Theorem 1]{FMS1})}.
\label{thm:FKE1-P}
  For a $Q$-function $q$ satisfying  condition \eqref{A-4} and for an arbitrary initial state distribution $\gamma\in \P(X),$  the jump process $\{\xi_t: t \in \BB{R}_+^0\}$ defined by  \eqref{JMP} on the probability space $(\Omega, \cal{F}, \F, \{\F_t\}_{t \ge 0}, \BB{P}_\gamma)$  is a jump Markov process with the transition function $\bar{P}.$
\end{theorem}

For simplicity, we write $\bar{P}(t,Z)$ instead of $\bar{P}(0,x,t,Z)$ when the initial state $x$ is fixed.
Let $E$ be a set and $\cal{A}$ be a set of functions $f:E\mapsto \bar{\BB{R}}.$ 
We say that $f$ is a minimal element of $\cal{A}$ if $f\in\cal{A}$ and $f(e)\le g(e)$ for all $e\in E$ and for all $g\in \cal{A}.$

Observe that, if
\begin{equation}
\label{boundedness}
\sup_{t \in \BB{R}_+} q(z,t) < +\infty, \qquad z \in X,
\end{equation}
then \eqref{A-4} holds. 
For a given $Q$-function $q,$  define the set of measurable subsets of $X$
\begin{equation}
\label{y}
\mathcal{Y} = \{Z \in \B(X): \sup_{z \in Z, t \in \BB{R}_+} q(z, t) < +\infty\}.
\end{equation}
The sets from $\mathcal{Y}$ are called  $q$-bounded; \cite{FMS1}. 
\begin{theorem} {\rm (\cite[Lemma~1(a) and Corollary~7 ]{FMS1}).}
\label{thm:FKE-P}
Fix an arbitrary $x \in X$. If \eqref{boundedness} holds, then

(a) there exists Borel subsets $X_n,$ $n = 1,2,\ldots,$ of $X$ such that $X_n \in \mathcal{Y}$ for all $n = 1,2,\ldots$ and $ X_n \uparrow X$ as $n \to \infty,$

(b) for all $t \in \BB{R}_+$ and $Z \in \mathcal{Y}$, the function $\bar{P}(t,Z)$ satisfies Kolmogorov's forward equation,
\begin{equation}
\label{MKDE-P}
P(t, Z) = \delta_{x}(Z) + \int_0^t  \int_X q(z,  s, Z\setminus \{z\}) P(s , dz)ds - \int_0^t \int_Z q(z, s)P(s,dz)ds.
\end{equation}

(c) the function $\bar{P}(t,Z)$, where $ t\in \BB{R}_+$ and $Z \in \B(X)$,
is the minimal function such that: (i) it is a measure on $(X, \B(X))$ for a  fixed $t$, (ii) it is measurable in $t$ for a fixed $Z$, and (iii) statement
(b) holds.  If $\bar{P}(s,X) = 1$ for some $s \in \BB{R}_+$, then $P(t,Z) = \bar{P}(t,Z),$ where $t \in ]0,s]$ and $Z \in \B(X),$
for every function $P:\BB{R}_+\times\B(X)\to [0,1]$ satisfying  conditions (i)--(iii). In addition, if $\bar{P}(t,X) =1 $ for all $t \in \BB{R}_+$, then $\bar{P}(t,Z)$ is the unique function with values in $[0,1]$ and satisfying  conditions (i)--(iii).
\end{theorem}

\begin{definition}\label{def:sol}
For  given $x\in X$ and $Q$-function $q,$ a function $P:\BB{R}_+\times X\mapsto \BB{R}$ is called a solution of Kolmogorov's forward equation \eqref{MKDE-P}, if $P$ satisfies properties (i)--(ii) stated in Theorem~\ref{thm:FKE-P}(c) and \eqref{MKDE-P} holds for all $t \in \BB{R}_+$ and $Z \in \mathcal{Y}.$
\end{definition}

Theorem~\ref{thm:FKE-P} states that $\bar{P}(\cdot,\cdot)$ is the minimal solution of Kolmogorov's forward equation \eqref{MKDE-P}.

Let $\phi$ be a Markov policy for a CTJMDP.  Consider the function $q(z,t,Z): = \tilde{q}(z, \phi_t, Z)$ defined by \eqref{q-def21} with $p(\cdot)= \phi(\cdot|z,t)$.  
Observe that the function $\tilde{q}(z,\phi_t,Z)$ is a $Q$-function, and, in view of Assumption~\ref{A1}, condition \eqref{boundedness} holds with $q(z,t) = \tilde{q}(z,\phi_t)$. In addition, the right-hand sides of formulae \eqref{m2} and \eqref{m22} coincide.  Therefore, $\nu^\phi=\nu.$  Therefore, $\BB{P}_\gamma^\phi= \BB{P}_\gamma$ for every initial state distribution $\gamma\in \P(X).$


 Let ${P}^\phi(u,x;t,Z)$, where  $u \in \BB{R}_+^0$, $x \in X$, $t \in ]u,+\infty[$, and $Z \in \B( X)$, be the transition function obtained from \eqref{b0}-\eqref{edef} with $q(z,t,Z) = \tilde{q}(z,\phi_t,Z)$. The following corollary from Theorem~\ref{thm:FKE-P} states that the jump process corresponding to the Markov policy $\phi$ is a jump Markov process with the transition function ${P}^\phi$. 

\begin{corollary}
\label{thm:FKE1}
For an initial distribution $\gamma$ on $X$ and a Markov policy $\phi$, the jump process $\{\xi_t: t\in \BB{R}_+^0\}$ defined on the probability space  $(\Omega, \cal{F}, \F, \{\F_t\}_{t \ge 0}, \BB{P}^\phi_\gamma)$  is a jump Markov process with the transition function ${P}^\phi$.
\end{corollary}
\proof{Proof.}
Recall that \eqref{boundedness} holds with $q(z,t) = \tilde{q}(z,\phi_t)$. Therefore,  inequality \eqref{A-4} holds with $q(z,s)=\tilde{q}(z,\phi_s)$, and the corollary follows from Theorem~\ref{thm:FKE1-P} with $q(z,t,Z) = \tilde{q}(z,\phi_t,Z),$ $\bar{P}(t,Z)=P^\phi(t,Z),$ and $\BB{P}_\gamma=\BB{P}_\gamma^\phi.$
\hfill \Halmos\endproof

Corollary~\ref{thm:FKE1} implies that, for the given Markov policy $\phi,$ initial state $x\in X,$ time epoch $t\in \BB{R},$ and set $Z\in\B(X),$
\begin{equation}\label{eq:ppx}
  P^\phi_x(t,Z)=\BB{P}_x^\phi(\xi_t\in Z)=P^\phi(0,x;t,Z),
\end{equation}
where the first equality follows from the definition of $P^\phi_x(t,Z)$ in \eqref{mm} and the second equality holds because $P^\phi$ is the transition function of the jump Markov process $\xi_t.$ 

Consider Kolmogorov's forward equation~\eqref{MKDE-P} with $q(z,t,Z) = \tilde{q}(z, \phi_t, Z)$,
\begin{equation}
\label{MKDE}
P(t, Z) = \delta_{x}(Z) + \int_0^t  \int_X \tilde{q}(z,  \phi_s, Z\setminus \{z\}) P(s , dz)ds - \int_0^t \int_Z \tilde{q}(z, \phi_s)P(s,dz)ds.
\end{equation}
Let $\mathcal{Y}^\phi:=\mathcal{Y},$ where $\mathcal{Y}$  is defined by \eqref{y} for $q(z,t)= \tilde{q}(z,\phi_t)$,
\begin{equation}
\label{y-varphi}
\mathcal{Y}^\phi = \{Z \in \B(X): \sup_{z \in Z, t \in \BB{R}_+} \tilde{q}(z, \phi_t) < +\infty\}.
\end{equation}
%
%

The following corollary from Theorem~\ref{thm:FKE-P} states that $P^\phi(\cdot,\cdot)$ is the minimal solution of Kolmogorov's forward equation~\eqref{MKDE}.

\begin{corollary}
\label{thm:FKE}
For fixed initial state $x \in X$ and  Markov policy $\phi$, the function $P_x^\phi(t,Z),$ defined for $t \in \BB{R}_+$ and $Z \in \B(X),$ is the minimal solution of Kolmogorov's forward equation~\eqref{MKDE}. If $P_x^\phi(s,X) = 1$ for some $s \in \BB{R}_+$, then $P(t,Z) = P_x^\phi(t,Z),$ where $t \in ]0,s]$ and $Z \in \B(X),$
for every solution $P(\cdot,\cdot)$ of Kolmogorov's forward equation~\eqref{MKDE} with values in $[0,1].$
In particular, if $P_x^\phi(t,X) = 1$ for all $t \in \BB{R}_+$, then $P_x^\phi(t,Z)$  is the unique solution of Kolmogorov's forward equation~\eqref{MKDE} with values in $[0,1]$.
\end{corollary}
\proof{Proof.} The corollary follows from Theorem~\ref{thm:FKE-P} with $q(z,t,Z) = \tilde{q}(z, \phi_t, Z),$ for all $z\in X,$ $t\in\BB{R}_+,$ and $Z\in\B(X),$ since $\bar{P}(t,Z)=\bar{P}(0,x,t,Z)=P^\phi(0,x,t,Z)=P_x^\phi(t,Z)$ for the $Q$-function $q,$ where the last equality follows from \eqref{eq:ppx}. \hfill\Halmos\endproof

\begin{corollary}
\label{C:FKE}
Let the transition rates $q(z,a)$ be bounded in $(z,a) \in Gr(A)$.
Then, for an initial state $x \in X$ and a Markov policy $\phi$, the function $P_x^\phi(t,Z),$ where $t \in \BB{R}_+$ and $Z \in \B(X)$, is the unique solution of Kolmogorov's forward equation~\eqref{MKDE}.
\end{corollary}
\proof{Proof.}
In view of \eqref{nonexp},  $P_x^\phi(t,X) = 1$ for all $t \in \BB{R}_+$ as the transition rates $q(z,a)$ are bounded. Thus, the corollary follows from the last conclusion of Corollary~\ref{thm:FKE}.
\hfill \Halmos\endproof

\section{Proof of Theorem~\ref{main}.}
\label{S-PM}

The proof of Theorem~\ref{main} is provided in two steps. First, we prove Lemma~\ref{L:i-s} stating that Theorem~\ref{main} holds when the initial distribution $\gamma$ on $X$ is concentrated at a point $x \in X.$ Second, using this fact, we prove that the theorem holds for all initial distributions $\gamma$ on $X$. The proof of Lemma~\ref{L:i-s} is based on the following lemma  stating that, for an initial state $x$ and for policies $\pi$ and $\varphi$ described in the statement of Theorem~\ref{main}, the marginal distributions $P_x^\pi(t,Z)$ and $P_x^\varphi(t,Z)$ are solutions to  Kolmogorov's forward equation~\eqref{MKDE} for the Markov policy $\phi=\varphi$.

\begin{lemma}
\label{solutions}
For an initial state $x \in X$ and for a policy $\pi$, let $\varphi$ be a Markov policy satisfying \eqref{formula-de2018} with $\gamma(\{x\}) =1$.  Then, the following statements hold:


(i) the functions $P_x^\pi(t, Z)$ and $P_x^\varphi(t,Z)$ are solutions of Kolmogorov's forward equation~\eqref{MKDE} with $\phi = \varphi$;

 (ii) for all $t \in \BB{R}_+$ and $Z \in \B(X)$,
 \begin{equation}
 \label{1-d}
  P_x^\varphi(t, Z) \le P_x^\pi(t, Z);
 \end{equation}

 (iii) if $P_x^\varphi(s,X) = 1$ for some $s \in \BB{R}_+$, then \eqref{1-d} holds for $t \in ]0,s]$ with an equality. In addition, if $P_x^\varphi(t,X) = 1$ for all $t \in \BB{R}_+$, then inequality~\eqref{1-d} holds with an equality for all $t \in \BB{R}_+$.
\end{lemma}
The proof of Lemma~\ref{solutions} is provided after presenting auxiliary Lemmas~\ref{l:Kitaev} and \ref{l:GKE}. Lemma~\ref{l:Kitaev} was introduced in Kitaev~\cite[Proof of Lemma~4]{Kit} and explicitly stated in Kitaev and Rykov~\cite[Lemma 4.28]{KR}. 
We provide the formulation of Lemma~\ref{l:Kitaev} here for completeness.

\begin{lemma}
\label{l:Kitaev} {\rm (Kitaev and Rykov~\cite[Lemma 4.28]{KR})}.
For an initial state $x \in X$ and for a policy $\pi$, consider the probability space $(\Omega, \F, \{\F_t\}_{t \ge 0}, \BB{P}_x^\pi),$  whose elements are defined in Section~\ref{S-MD}.
Then the random measure 
\begin{equation}
\label{m4}
\tilde{\nu}^\pi( [0,t], Z) := \int_0^t \tilde{q}(\xi_{s}, \pi_s) I\{\xi_{s}\in Z\} ds, \qquad t \in \BB{R}_+^0,\ Z \in \B(X),
\end{equation}
is a compensator of the random measure 
\begin{equation}
\label{m3}
\tilde{\mu}( [0,t],Z) := \sum_{n \geq 1} I\{ t_n \in\, [0,t]\}I\{x_{n-1} \in Z\}, \qquad t \in \BB{R}_+^0,\ Z \in \B(X).
\end{equation}
\end{lemma}
We remark that $\tilde{\mu}( [0,t],Z) $ is the number of jumps out of the set $Z$  and $\tilde{\nu}^\pi( [0,t], Z)$ is the cumulative intensity of jumping out of the set $Z$ during the time interval $]0,t]$.  We recall that the random measures ${\mu}( [0,t],Z) $ and ${\nu}^\pi( [0,t], Z)$ described in \eqref{m1} and \eqref{m2} deal with the numbers of jumps to sets $Z\in\B(X).$

\begin{definition}
\label{D:xpi}
For an initial state $x \in X$ and a policy $\pi$, a set $Z \in \B(X)$ is called $(x, \pi)$-bounded if $\sup_{t\in \BB{R}_+} \BB{E}_x^\pi \tilde{q}(\xi_t, \pi_t) I\{\xi_t \in Z\} < +\infty.$
\end{definition}

The following lemma  implies that,  for an initial state $x\in X$ and a policy $\pi$, the marginal distribution  $P_x^\pi(t,Z)$ satisfies, for all $t \in \BB{R}_+$ and $(x, \pi)$-bounded sets $Z \in \B(X)$, the equation
\begin{equation}
\label{GKDE}
P_x^\pi(t, Z) = \delta_{x}(Z)+ \BB{E}_x^\pi \int_{0}^{t} \tilde{q}(\xi_{s}, \pi_s, Z\setminus \{\xi_s\})I\{\xi_s \in X\} ds - \BB{E}_x^\pi \int_{0}^{t} \tilde{q}(\xi_{s}, \pi_s)I\{\xi_s \in Z\} ds.
\end{equation}
For bounded transition rates,  Kitaev~\cite[Lemma 4]{Kit} showed that \eqref{GKDE} holds for all $Z \in \B(X).$ 
For unbounded transition rates, the integrals in \eqref{GKDE} need not be finite, and hence their difference may not be defined. However, if  $\sup_{z \in Z} \bar{q}(z) < +\infty$ for $Z \in \B(X)$, then
\begin{equation}
\label{xpi-q}
\BB{E}_x^\pi \tilde{q}(\xi_{s}, \pi_s)I\{\xi_s \in Z\}  \le \BB{E}_x^\pi \bar{q}(\xi_s) I\{\xi_s \in Z\} \le \sup_{z \in Z}\bar{q}(z) < +\infty, \quad t \in \BB{R}_+,
\end{equation}
where the first inequality follows from \eqref{f1} with $p(\cdot) = \pi_s(\cdot)$, and the other inequalities are straightforward. Therefore, the second integral in \eqref{GKDE} is finite for sets $Z \in \B(X)$ with $\sup_{z \in Z} \bar{q}(z) < +\infty$. Using this fact, Guo and Song~\cite[Theorem 3.1(c)]{GS} and Piunovskiy and Zhang~\cite[Theorem 1(b)]{PZ1} showed that \eqref{GKDE} holds for  sets $Z \in \B(X)$ such that $\sup_{z \in Z} \bar{q}(z) < +\infty$ when the transition rates are unbounded and the associated jump process is nonexplosive. The condition for nonexplosiveness considered in \cite[Condition 1]{PZ1}  is more general than the condition considered in \cite[Assumption A]{GS}. In the following lemma, we show that formula~\eqref{GKDE} holds for possibly explosive jump processes and for $(x,\pi)$-bounded sets, which is a larger class of sets than the sets considered in \cite{PZ1}.  Therefore, Lemma~\ref{l:GKE} generalizes \cite[Theorem 1(b)]{PZ1}; see Corollary~\ref{c:GKE}.

\begin{lemma}
\label{l:GKE}
For an initial state $x \in X$ and a policy $\pi$, 
formula~\eqref{GKDE} holds for all $t \in \BB{R}_+$ if the set $Z\in \B(X)$ is $(x, \pi)$-bounded.
\end{lemma}
\proof{Proof.}
For all $m = 0,1, \ldots$, $t \in \BB{R}_+$, and $Z \in \B(X)$,  the number of jumps in the interval $[0,t\wedge t_m]$ is bounded by $m$. Then, as the random measures $\mu([0,t\wedge t_m],Z)$ and $\tilde{\mu}([0,t\wedge t_m], Z)$ defined in \eqref{m1} and \eqref{m3}, respectively, give the numbers of jumps into and out of the set  $Z$ during the interval $]0,t\wedge t_m]$,
\begin{equation}
\label{PEF}
I\{\xi_{t\wedge t_m}\in Z\} = \delta_x(Z) + \mu([0,t\wedge t_m],Z) - \tilde{\mu}([0,t\wedge t_m], Z) .
\end{equation}
Taking expectation in the both sides of \eqref{PEF} implies that, for all $m = 0,1, \ldots$, $t \in \BB{R}_+$, and $Z \in \B(X)$,
\begin{equation}
\label{P-E}
\BB{P}_x^\pi(\xi_{t\wedge t_m}\in Z) = \delta_{x}(Z) + \BB{E}_x^\pi (\mu([0,t\wedge t_m],Z)) - \BB{E}_x^\pi (\tilde{\mu}([0,t\wedge t_m], Z)).
\end{equation}

To prove \eqref{GKDE} for $(x, \pi)$-bounded sets $Z \in \B(X)$, we take $m \to \infty$ in formula~\eqref{P-E}. For all $t \in \BB{R}_+$ and $Z \in \B(X)$,
\begin{equation}
\label{inc1}
\lim_{m \to \infty} \BB{E}_x^\pi (\mu([0,t\wedge t_m], Z) ) = \BB{E}_x^\pi (\mu([0,t\wedge t_\infty], Z)) = \BB{E}_x^\pi (\mu([0,t], Z)) = \BB{E}_x^\pi (\nu^\pi([0,t],Z)),
\end{equation}
where the first equality follows from the monotone convergence theorem since $\mu([0,t\wedge t_m], Z) \uparrow \mu([0,t\wedge t_\infty], Z)$ as $m \to \infty$, the second equality is correct because $\{t_{n} \in [0,t]\} = \{t_n \in [0,t\wedge t_\infty]\}$ for all $n = 1,2,\ldots$, and the last one follows from \eqref{mu-nu} with $\BB{E} = \BB{E}_x^\pi$, $\nu = \nu^\pi$, and $T = t$  since $\nu^\pi$ is a compensator of the random measure $\mu.$ Similarly, since $\tilde{\mu}([0,t\wedge t_m], Z)  \uparrow \tilde{\mu}([0,t\wedge t_\infty], Z)$ as $m \to \infty$ and, in view of  Lemma~\ref{l:Kitaev},  $\tilde{\nu}^\pi$ is a compensator of the random measure $\tilde{\mu}$,
\begin{equation}
\label{inc2}
\lim_{m \to \infty} \BB{E}_x^\pi (\tilde{\mu}([0,t\wedge t_m], Z) ) = \BB{E}_x^\pi (\tilde{\mu}([0,t\wedge t_\infty], Z)) = \BB{E}_x^\pi (\tilde{\mu}([0,t], Z)) = \BB{E}_x^\pi (\tilde{\nu}^\pi([0,t],Z)).
\end{equation}

Let us fix an arbitrary $t\in \BB{R}_+,$ a policy $\pi,$ and an $(x,\pi)$-bounded set $Z \in \B(X).$ Observe that
\begin{equation}
\label{t-u-bound}
\BB{E}_x^\pi (\tilde{\nu}^\pi([0,t], Z)) =  \int_{0}^{t}   \BB{E}_x^\pi \tilde{q}(\xi_{s}, \pi_s) I \{ \xi_{s} \in Z\}ds
\le t\left (\sup_{s \in \BB{R}_+} \BB{E}_x^\pi \tilde{q}(\xi_{s}, \pi_s) I \{ \xi_{s} \in Z\} \right) < +\infty,
\end{equation}
where the first equality follows from \eqref{m4} and from interchanging the expectation and integration operators, 
the  first inequality is straightforward, and the last one holds since $Z$ is an $(x,\pi)$-bounded set.
Then
\begin{multline}
\label{rlimit}
\lim_{m \to \infty} (\BB{E}_x^\pi (\mu([0,t\wedge t_m],Z)) - \BB{E}_x^\pi (\tilde{\mu}([0,t\wedge t_m], Z)))\\
\begin{aligned}
&= \BB{E}_x^\pi (\nu^\pi([0,t] \times Z))-  \BB{E}_x^\pi (\tilde{\nu}^\pi([0,t] \times Z)) \\
&=  \BB{E}_x^\pi \int_{0}^{t}\tilde{q}(\xi_{s}, \pi_s, Z \setminus \{\xi_{s}\} )I\{\xi_s \in X\}ds - \BB{E}_x^\pi  \int_{0}^{t}\tilde{q}(\xi_{s}, \pi_s) I \{ \xi_{s} \in Z\}ds,
\end{aligned}
\end{multline}
where, in view of \eqref{t-u-bound}, the first equality follows from \eqref{inc1} and \eqref{inc2}, and
the last one follows from (\ref{m2}) and (\ref{m4}).

In addition, 
\begin{multline}
\label{llimit}
\begin{aligned}
\lim_{m \to \infty}\BB{P}_x^\pi(\xi_{t \wedge t_m} \in Z) &= \lim_{m \to \infty} (\BB{P}_x^\pi(\xi_{t} \in Z, t < t_m) + \BB{P}_x^\pi(\xi_{t_m} \in Z, t \ge t_m)) \\
&= \BB{P}_x^\pi(\xi_{t} \in Z, t < t_\infty) = \BB{P}_x^\pi(\xi_t \in Z) = P_x^\pi(t,Z),
\end{aligned}
\end{multline}
where the first equality is straightforward,  the third equality holds because $\{\xi_t \in X\} = \{t < t_\infty\}$, the last one is correct due to \eqref{mm}, and the second equality holds in view of the continuity  of probability because $\{\xi_{t} \in Z, t < t_m\}\downarrow \{\xi_{t} \in Z, t < t_\infty\}$ as $m \to \infty$ and, as shown in the rest of this proof,
\begin{equation}
\label{l-0}
\lim_{m \to \infty}\BB{P}_x^\pi(\xi_{t_m} \in Z, t \ge t_m) = 0
\end{equation}
 for an  $(x,\pi)$-bounded set $Z.$  To complete the proof of \eqref{llimit}, we need to verify \eqref{l-0}.

Let  $\limsup_{m \to \infty}\BB{P}_x^\pi(\xi_{t_m} \in Z, t \ge t_m)= p > 0$. Then there exists a subsequence $\{m_k, k = 1,2,\ldots\}$ such that $\BB{P}_x^\pi(\xi_{t_{m_k}}\in Z, t \ge t_{m_k}) > \frac{p}{2}$ for all $k = 1,2,\ldots\ $. This fact and \eqref{m1} imply that
\begin{equation}
\label{contra1}
\BB{E}_x^\pi (\mu([0,t], Z)) = \sum\limits_{m \ge 1} \BB{P}_x^\pi (t_{m} \in [0,t], x_{m} \in Z) \ge \sum\limits_{k \ge 1} \BB{P}_x^\pi (t_{m_k} \in [0,t], x_{m_k} \in Z)= +\infty.
\end{equation}
Since the set $Z$ is  $(x,\pi)$-bounded,
\begin{equation}\label{eqendogprKE} \BB{E}_x^\pi (\mu([0,t], Z])) = \lim_{m\to \infty} \BB{E}_x^\pi(\mu([0,t\wedge t_m], Z)) \le \lim_{m \to \infty} \BB{E}_x^\pi(\tilde{\mu}([0,t\wedge t_m], Z)) + 1 < +\infty,\end{equation}
where the first equality follows from the first and second equalities in \eqref{inc1}, the first inequality follows from \eqref{P-E}, and the last inequality follows from \eqref{inc2} and \eqref{t-u-bound}. Since inequality~\eqref{eqendogprKE} contradicts \eqref{contra1}, formula~\eqref{l-0} holds, and \eqref{llimit} is proved.

Let $m\to\infty$ in \eqref{P-E}. In view of \eqref{llimit}, the left-hand side of \eqref{P-E} tends to $P_x^\pi(Z).$ In view of \eqref{rlimit}, the right-hand side of \eqref{P-E} tends to the right-hand side of \eqref{GKDE}.  So, \eqref{GKDE} is proved.
 \hfill\Halmos\endproof

The following corollary generalizes Piunovskiy and Zhang~\cite[Theorem 1(b)]{PZ1} since it holds for possibly explosive jump processes. 
\begin{corollary}{\rm (cp. Piunovskiy and Zhang~\cite[Theorem 1(b)]{PZ1})}
\label{c:GKE}
For an initial state $x \in X$ and a policy $\pi$, formula~\eqref{GKDE} holds for all $t \in \BB{R}_+$ and for all $Z \in \B(X)$ with $\sup_{z \in Z}\bar{q}(z) < +\infty$.
\end{corollary}
\proof{Proof.}
As follows from \eqref{xpi-q}, every set $Z \in \B(X)$ with $\sup_{z \in Z}\bar{q}(z) < +\infty$ is $(x,\pi)$-bounded. Thus, the corollary follows from Lemma~\ref{l:GKE}.
 \hfill\Halmos\endproof

\proof{Proof of Lemma~\ref{solutions}.}
Consider the set of measurable sets $\mathcal{Y}^\varphi$ defined by \eqref{y-varphi}. In view of Corollary~\ref{thm:FKE}, statements (ii) and (iii) of the lemma hold if statement (i)  holds. The rest of the proof verifies statement (i) of the lemma.

Corollary~\ref{thm:FKE} with $\phi = \varphi$ implies that the function $P_x^\varphi(t,Z)$ is the minimal solution of Kolmogorov's forward equation~\eqref{MKDE}. It remains to show that the function $P_x^\pi(t,Z)$ is also a solution of \eqref{MKDE}. Observe that, for any non-negative measurable function $f$, for all $Z \in \B(X)$ and $s \in \BB{R}_+$,
\begin{multline}
\label{f-int0}
\begin{aligned}
\BB{E}_x^\pi f(\xi_s, \pi_s) I\{\xi_s \in Z\}  &= \int_Z \int_{A(z)} f(z,a) P_x^\pi(s,dz,da)\\
&= \int_Z \int_{A(z)} f(z,a) \varphi(da|z,s) P_x^\pi(s,dz) = \int_Z  f(z,\varphi_s)  P_x^\pi(s,dz),
\end{aligned}
\end{multline}
where the first equality follows from \eqref{f-Ex}, the second equality follows from \eqref{formula-de2018}, and the last one follows from \eqref{Ext}. Then, for any non-negative measurable function $f$, for all $t \in \BB{R}_+$ and $Z \in \B(X)$,
\begin{equation}
\label{f-int}
\BB{E}_x^\pi \left(\int_0^t f(\xi_s, \pi_s) I\{\xi_s \in Z\} ds\right) =  \int_0^t  \BB{E}_x^\pi f(\xi_s, \pi_s) I\{\xi_s \in Z\} ds = \int_0^t \int_Z  f(z,\varphi_s)  P_x^\pi(s,dz) ds,
\end{equation}
where the first equality  follows from interchanging integration and expectation, and the second one follows from \eqref{f-int0}. Therefore, it follows from Lemma~\ref{l:GKE}, formula \eqref{f-int} with $Z = X$ and $f(\xi_s, \pi_s) = \tilde{q}(\xi_s, \pi_s, Z \setminus \{\xi_s\})$, and the same formula with $f(\xi_s, \pi_s) = \tilde{q}(\xi_s, \pi_s)$ imply that the function $P_x^\pi(t,Z)$ satisfies Kolmogorov's forward equation~\eqref{MKDE} with $\phi = \varphi$ for all $t \in \BB{R}_+$ if the set $Z \in \B(X)$ is  $(x,\pi)$-bounded.

To conclude the proof of the lemma, observe that every set $Z \in \mathcal{Y}^\varphi$ is $(x,\pi)$-bounded, where $\mathcal{Y}^\varphi$ is defined in \eqref{y-varphi} for $\phi=\varphi.$ This is true because, for $Z \in \mathcal{Y}^\varphi$,
\begin{multline*}
\label{bound}
\sup_{s \in \BB{R}_+}\BB{E}_x^\pi  \tilde{q}(\xi_s, \pi_s) I\{\xi_s \in Z\} =  \sup_{s \in \BB{R}_+} \int_Z \tilde{q}(z,\varphi_s) P_x^\pi(s,dz) \le \left(\sup_{z \in Z, s \in \BB{R}_+} \tilde{q}(z,\varphi_s) \right)P_x^\pi (s, Z) <  +\infty,
\end{multline*}
where the first equality follows from \eqref{f-int0} with $f(\xi_s, \pi_s) = \tilde{q}(\xi_s, \pi_s)$, the first inequality is straightforward, and the last one is true since $Z \in \mathcal{Y}^\varphi$ and $P_x^\pi(s, Z) \le 1$. Therefore, the function $P_x^\pi(t,Z)$ is a solution of Kolmogorov's forward equation~\eqref{MKDE}, and statement (i) of the lemma holds.
\hfill \Halmos\endproof

The following lemma is Theorem~\ref{main} with the initial distribution $\gamma$ concentrated at a point $x \in X$.
\begin{lemma}
\label{L:i-s}
Theorem~\ref{main} holds if the initial distribution $\gamma\in\P(X)$ is concentrated at a singleton, that is,  $\gamma(\{x\})= 1$ for some $x \in X.$
\end{lemma}
\proof{Proof.}
For all $t \in \BB{R}_+$, $Z \in \B(X) $, and $B \in \B(A)$,
\begin{equation}
\label{2-d-e}
\begin{aligned}
P_x^\varphi(t, Z, B) &= \BB{E}_x^\varphi \left[ I\{\xi_{t-} \in Z\} \varphi_t(B) \right]  
  =  \int_Z \varphi(B | z, t) P_x^\varphi(t, dz)\\
  & \le \int_Z \varphi(B | z, t) P_x^\pi(t, dz) = P_x^\pi(t, Z, B),
\end{aligned}
\end{equation}
where the first equality is \eqref{M} with $\pi=\varphi$, the second equality follows from \eqref{mm}, \eqref{zero-m}, and the property  $\varphi(\,\cdot\, | \omega, t) = \varphi(\,\cdot\, | \xi_{t-}, t)$  for the Markov policy $\varphi$,  the inequality is correct since \eqref{1-d} holds as follows from Lemma~\ref{solutions}(ii), and the last equality is \eqref{formula-de2018}. Therefore, \eqref{2-d} holds.

Observe that \eqref{2-d-e}  holds with an equality if $P_x^\varphi(t,Z) = P_x^\pi(t,Z)$ for all $Z \in \B(X)$. Therefore, in view of Lemma~\ref{solutions}(iii), formula \eqref{2-d} holds with an equality for $t\in ]0,s],$ if $P_x^\varphi(s,X) = 1$ for some $s \in \BB{R}_+$, and for $t \in \BB{R}_+$ if $P_x^\varphi(t,X) = 1$ for all $t \in \BB{R}_+$.
 \hfill\Halmos\endproof

\begin{definition}
A sub-model of a CTJMDP $\{X',A',A'(\cdot),q'\}$ is a CTJMDP $\{X,A,A(\cdot),\tilde{q}\}$ with $X \subseteq X'$, $A \subseteq A'$, $A(z) \subseteq A'(z)$ for all $z \in X$, and $\tilde{q}(z,a,Z) = q'(z,a,Z)$ for all $z \in X$, $a \in A(z),$ and $Z \in \B(X)$.
\end{definition}

Theorem~\ref{main} is proved in Lemma~\ref{L:i-s} for the initial distribution concentrated at one point.  The extension to an arbitrary initial distribution is based on the following arguments explained in detail in the proof of Theorem~\ref{main}.  We introduce an extended CTJMDP, for which our original CTJMDP is a sub-model, by adding a state $x'$ to $X$ and actions $a',a''$ to $A$. 
Then, for an arbitrary policy $\sigma$ for the original CTJMDP and for an arbitrary fixed $u \in \BB{R}_+,$ we construct in  a natural way a policy $\tilde{\sigma}$ for the extended CTJMDP such that  the marginal distributions satisfy
 \begin{equation}\label{eq5.16EF28}
 P_{x'}^{\tilde{\sigma}}(t+u, Z,B) = (1-e^{-u})P_\gamma^\sigma(t,Z,B), \qquad t \in \BB{R}_+, Z \in \B(X), B \in \B(A).
 \end{equation}

For an arbitrary policy $\pi$ and the corresponding Markov policy $\varphi$ satisfying \eqref{formula-de2018}, we shall consider policies $\tilde{\pi}$ and $\tilde{\varphi}$ in the extended CTJMDP such that formula \eqref{eq5.16EF28} holds for $(\sigma, \tilde{\sigma})= (\pi, \tilde{\pi})$, and for $(\sigma,\tilde{\sigma})=(\varphi,\tilde{\varphi})$. We shall also show that for the extended CTJMDP formula \eqref{formula-de2018} holds being applied to  the policy $\tilde{\pi}$, Markov policy $\tilde{\varphi}$, and initial distribution concentrated at the state $x'.$ Then, as follows from  Lemma~\ref{L:i-s} applied to the extended CTJMDP,
\begin{equation}
\label{2-ds}
P_{x'}^{\tilde{\varphi}} (t, Z, B) \le P_{x'}^{\tilde{\pi}} (t, Z, B), \qquad t \in \BB{R}_+, Z \in \B(X), B \in \B(A).
\end{equation}
In view of this fact, formula \eqref{eq5.16EF28}, applied to the pairs of policies $(\sigma,\tilde{\sigma})=(\pi,\tilde{\pi})$ and $(\sigma,\tilde{\sigma}) = (\varphi,\tilde{\varphi}),$ implies the correctness of Theorem~\ref{main}; see the diagram in Figure~\ref{figure}.

\begin{figure}[h]
\begin{equation*}
\begin{array}{ccc}
P_x^{\tilde{\varphi}}(t+u, Z, B) &\le &P_x^{\tilde{\pi}}(t+u, Z, B) \\
\verteq & &\verteq  \\
(1-e^{-u})P_\gamma^\varphi(t,Z,B) &  &(1-e^{-u})P_\gamma^\pi(t,Z,B)
\end{array}
\implies P_\gamma^\varphi(t, Z, B) \le P_\gamma^\pi(t, Z, B)
\end{equation*}
\caption{Major steps of the proof of Theorem~\ref{main}.}
\label{figure}
\end{figure}

\proof{Proof of Theorem~\ref{main}.}
Let us fix an arbitrary $u\in\BB{R}_+.$
For $x' \notin X$ and $a', a'' \notin A$,  let $X' := X \cup \{x'\}$, $A' := A \cup \{a', a''\}$, $A'(x) := A(x) \cup \{a''\}$ for all $x \in X,$ and $A'(x') :=\{ a', a''\}$.
For all $x \in X', a \in A'(x),$ and $Z \in \B(X')$, define the new transition rate $q'$ by
\begin{equation}
\label{new-q1}
q'( x,a,Z) := \left \{
\begin{array}{ll}
\tilde{q}( x, a, Z \setminus \{x'\} ), &\quad \text{ if } x \in X, a \in A(x),\\
\gamma(Z\setminus \{x'\})-\delta_{x'}(Z)   &\quad \text{ if }  x = x', a = a',\\
0, &\quad \text{ if } x \in X', a = a''.
\end{array}
\right.
\end{equation} This means that an additional state $x'$ is added to the state space $X.$  The set of feasible actions at this state consists of  two actions $a'$ and $a''.$ If the action $a'$ is chosen at the state $x',$ then the process jumps to every set $Z\in\B(X)$ with the intensity $\gamma(Z).$ In addition, the  action $a''$ is also added to the action sets $A(x)$ for all $x\in X.$ Under this action, every set $x\in X'=X\cup\{x'\}$ is absorbing.

Consider  the extended CTJMDP $\{X', A', A'(x), q'\}$. For this extended CTJMDP,  we shall respectively denote by $x_n',$  $t_{n}',$ $\omega'$, $\Omega'$, $\xi_t',$ and $U_t'$  the objects $x_n,$  $t_{n}$,  $\Omega$, $\xi_t,$ and $U_t$  defined in Section \ref{S-MD},   where $n = 0,1,\ldots.$ In particular, consider the set of trajectories ${\Omega}'(u)$ defined by
\[{\Omega}'(u) = \{(x', t_1', x_0, t_1+u, x_1, t_2+u, \ldots):\, t_1' \in ]0,u],\  (x_0, t_1, x_1, t_2, \ldots) \in \Omega\}.\]
For $\omega'=(x'_0,t'_1x'_1,t'_2,\ldots)\in \Omega'(u),$ let us define $\omega'_{-u}:=(x'_1,(t'_2-u),x'_2,(t'_3-u),\ldots)$, which is  the sample path starting from time $u$ and shifted back by time $u.$   Note that $\xi_{t}(\omega) = \xi_{t+u}'(\omega')$ for all $t \in \BB{R}_+,$ $\omega \in \Omega,$ and $\omega' \in {\Omega'}(u)$. The definition of ${\Omega}'(u)$ implies that it is a measurable subset of $\Omega'.$ 

For a policy $\sigma$ in the original CTJMDP, let $\tilde{\sigma}$ be a policy for the extended CTJMDP such that, for all $\omega' \in \Omega'$, $t \in \BB{R}_+$, and $B \in \B(A')$,
\begin{equation}
\label{new-p}
\tilde{\sigma}(B \mid \omega',t) = \left\{ \begin{array}{ll}
\delta_{a'}(B) I\{\xi_{t-}' = x'\} + \delta_{a''} (B)I\{\xi_{t-}' \in X\}, & \quad \text{ if } t \in ]0,u], \\
\delta_{a''} (B)I\{\xi_{t-}' = x'\}, & \quad \text{ if } t \in ]u, +\infty[,\\
\sigma(B \setminus \{a',a''\}\mid \omega'_{-u}, t-u)I\{\omega' \in \Omega'(u)\}, & \quad \text{ if } t \in ]u, +\infty[.
\end{array}
\right.
\end{equation}
 For the initial state  $x'_0=x',$  the policy $\tilde{\sigma}$  chooses the action $a'$ at the state $x'$ during the time interval $]0;u].$ If the jump does not occur during the time interval $]0;u],$ the policy $\tilde{\sigma}$  always chooses the action $a''$ at the state $x'$ during the time interval $]u,+\infty[.$  Of course, in this case  the state $x'$ becomes absorbing at the time instance $u$.  If the first jump occurs during the time interval $]0;u],$ then the process jumps to a state $x'_1\in X$.  Observe that, in view of \eqref{new-q1},
 \begin{equation}\label{eprelfff}
 \BB{P}_{x'}^{\tilde{\sigma} }(x'_1\in Z\mid\ t'_1\le u)=\gamma(Z) \qquad {\rm and}  \qquad \BB{P}_{x'}^{\tilde{\sigma}}(t'_1\le t)=1-e^{-t}
 \end{equation}
 for $Z \in \B(X)$ and for $t\in [0,u].$  After this jump the process stays at the state $x'_1\in X$ until the time epoch $u,$ and during the time interval $]u,+\infty)$ the policy $\tilde{\sigma}$ selects actions at time instances $t\in ]u,+\infty]$ in the same way   as the policy $\sigma$ does at time instances $t-u$ using the  observations starting  from the initial state $x_0=x'_1$ and initial time 0 until the time epoch $t-u.$  Thus, $\BB{P}_{x'}^{\tilde{\sigma}}(\xi'_{t+u}\in Z\mid t'_1\le u)=\BB{P}_{x'_1}^\sigma (\xi_t\in Z) $\ \   $\BB{P}_{x'}^{\tilde{\sigma}}$--a.s. for $t\in\BB{R}_+$ and for $Z\in\B(X).$  For                                                                                       $Z\in\B(X),$
   \begin{equation}\label{eqdescrx1}
 \BB{P}_{x'}^{\tilde{\sigma}}( x'_1\in Z)=\BB{P}_{x'}^{\tilde{\sigma}}(t'_1\le u,\, x'_1\in Z)=(1-e^{-u})\gamma(Z),
 \end{equation}
where the first equality holds because $\{t'_1\le u\}=\{ x'_1\in X\}$ up to null set of probability $\BB{P}_{x'}^{\tilde{\sigma}},$ and the second one follows from \eqref{eprelfff}.

For $t\in \BB{R}_+,$  $Z\in\B(X),$ and $B\in\B(A),$
 \begin{equation}\label{eqlohgEFFA}
 \begin{aligned}
& \quad P_{x'}^{\tilde{\sigma}}(t+u,Z,B)=\BB{E}_{x'}^{\tilde{\sigma}}[I\{x'_1\in X\}I\{\xi'_{t+u}\in Z\}\tilde{\sigma}(B\mid \omega',t+u)]\\
 &= \BB{E}_{x'}^{\tilde{\sigma}}[I\{x'_1\in X\}\BB{E}_{x'}^{\tilde{\sigma}}[I\{\xi'_{t+u}\in Z\}\tilde{\sigma}(B\mid \omega',t+u)|x'_1]]\\
 &=\BB{E}_{x'}^{\tilde{\sigma}}[I\{x'_1\in X\}\BB{E}_{x'_1}^\sigma [I\{\xi_{t}\in Z\}\sigma(B\mid \omega,t)]]\\
&=\BB{E}_{x'}^{\tilde{\sigma}}[I\{x'_1\in X\}P_{x'_1}^\sigma(t,Z,B)]=\int_X P_{x'_1}^\sigma(t,Z,B)\BB{P}_{x'}^{\tilde{\sigma}}(dx'_1)=(1-e^{-u})P_\gamma^\sigma(t,Z,B),
 \end{aligned}
 \end{equation}
 where  the first equality follows from  \eqref{M} and from $\{x'_1\in X\}=\{ \xi'_{t+u}\in X\}\supset \{ \xi'_{t+u}\in Z\}$ for $Z\in\B(X)$ and  the equality of sets holds up to null sets of measure  $\BB{P}_{x'}^{\tilde{\sigma}}.$
   The second equality in \eqref{eqlohgEFFA} follows from the properties of conditional expectations  because the function $I\{x'_1\in X\}$ is $\sigma(x'_1)$-measurable,  the third equality follows from the definition of the policy $\tilde{\sigma}$ in \eqref{new-p},  the fourth equality follows from \eqref{M},
   the fifth equation follows from the definitions of expectations and indicators, and the last one follows from \eqref{eqdescrx1}.

In view of \eqref{eqlohgEFFA}, for    $t\in \BB{R}_+,$  $Z\in\B(X),$ and $B\in\B(A),$
\begin{equation}
 \label{zero}
P_\gamma^\sigma(t,Z,B) = (1-e^{-u})^{-1}P_{x'}^{\tilde{\sigma}}(t+u,Z,B).
\end{equation}
Thus,  $P_\gamma^\sigma(t,\cdot)$ and $P_{x'}^{\tilde{\sigma}}(t+u,\cdot)$ are equivalent measures on $(X,\B(X))$ for all  $t\in \BB{R}_+.$

Let us fix an arbitrary policy $\pi$ for the original CTJMDP.  Let $\varphi$ be a Markov policy satisfying equality \eqref{formula-de2018}, and let us consider policies $\tilde{\pi}$ and $\tilde\varphi$ for the extended CTJMDP satisfying  \eqref{new-p} with $\sigma=\pi$  and $\sigma=\varphi$ respectively.  As follows from   \eqref{new-p} applied to $\sigma=\varphi,$ the Markov policy $\tilde\varphi$ is defined uniquely and,  for $B\in\B(A'),$ $z\in X',$ and $t\in\BB{R}_+,$
\begin{equation*}\begin{aligned}\tilde{\varphi}(B\mid z,t)=&\delta_{a'}(B)I\{z=x',t\le u\}+\delta_{a''}(B)I\{\{z\in X,t\le u\}\cup\{z=x',t>u\}\}\\ &+\varphi(B\setminus\{a',a''\}|z,t-u)I\{z\in X, t>u\}.\end{aligned} \end{equation*}
In view of \eqref{formula-d} and \eqref{zero} with $\sigma=\pi,$ for $B\in\B(A),$ $z\in X,$ and $t\in\BB{R}_+,$
\[
\tilde{\varphi}(B|z,t+u)=\varphi(B|z,t)=\frac{P_\gamma^\pi(t,dz,B)}{P_\gamma^\pi(t,dz)}= \frac{P_{x'}^{\tilde{\pi}}(t+u,dz,B)}{P_{x'}^{\tilde{\pi}}(t+u,dz)}, \qquad (P_{x'}^{\tilde{\pi}}(t+u,\cdot)\emph{-a.s.}), 
\]
which implies that formula \eqref{formula-d} folds for the policies $\tilde\varphi,$ $\tilde{\pi}$ and the initial state distribution concentrated at the state $x'$ when $z\in X$ and $t>u.$  Formula \eqref{formula-d} also holds for $\tilde\varphi,$ $\tilde{\pi},$  for the initial state distribution concentrated at the state $x',$  and for state-time pairs $(z,t)\in X\times\BB{R}_+\setminus\{z\in X,t>u\}$
 because at these state-time pairs the policy  $\tilde\pi$ is Markov. Since \eqref{formula-d} is equivalent to \eqref{formula-de2018},  for $t\in\BB{R}_+,$  $Z\in\B(X),$ and $B\in\B(A),$
\begin{equation}\label{elastcaseEF}
P_\gamma^\varphi(t,Z,B) = (1-e^{-u})^{-1}P_{x'}^{\tilde{\varphi}}(t+u,Z,B)\le (1-e^{-u})^{-1}P_{x'}^{\tilde{\pi}}(t+u,Z,B)=P_\gamma^\pi(t,Z,B),
\end{equation}
where the inequality follows from Lemma~\ref{L:i-s} and the equalities follow from \eqref{zero} applied to $\sigma=\varphi$ and $\sigma=\pi$ respectively. Inequality \eqref{2-d} is proved.

To complete the proof of the theorem, assume that $P_\gamma^\varphi(s,X) = 1$ for some $s \in \BB{R}_+$. We fix an arbitrary $t\in ]0,s].$ Then $P_\gamma^\varphi(t,X) = 1.$

Let $P_\gamma^\varphi(t,Z,B)<P_\gamma^\pi(t,Z,B)$ for some $Z\in\B(X)$ and $B\in\B(A).$ Then, in view of \eqref{elastcaseEF}, $P_\gamma^\varphi(t,Z,A\setminus B)\le P_\gamma^\pi(t,Z,A\setminus B)$ and
$P_\gamma^\varphi(t,X\setminus Z,A)\le P_\gamma^\pi(t,X\setminus Z,A).$  Therefore
\[\begin{aligned}
&1= P_\gamma^\varphi(t,X)=P_\gamma^\varphi(t,Z,B)+P_\gamma^\varphi(t,Z,A\setminus B)+P_\gamma^\varphi(t,X\setminus Z,A)\\ &<P_\gamma^\pi(t,Z,B)+P_\gamma^\pi(t,Z,A\setminus B)+P_\gamma^\pi(t,X\setminus Z,A)=P_\gamma^\pi(t,X),
\end{aligned}
\]
which is impossible since $P_\gamma^\pi(t,X)\le 1$.  Thus,   $P_\gamma^\varphi(t,Z,B)=P_\gamma^\pi(t,Z,B)$ for all $Z\in\B(X)$ and for all $B\in\B(A).$ \hfill\Halmos\endproof

\section{Sufficiency of Markov policies for  particular objective criteria.} \label{S-ER}
This section describes applications of the results of Section~\ref{S-MR} to CTJMDPs with finite and infinite horizons.  For finite-horizon CTJMDPs we consider expected total discounted costs.  For infinite-horizon CTJMDPs we consider expected total discounted costs and average costs per unit of time.
 For each of these problems we show that for a fixed arbitrary initial distribution $\gamma$  the objective criterion for a Markov policy $\varphi,$ described in Theorem~\ref{main} for a policy $\pi,$ is smaller than or equal to  the objective criterion for the policy $\pi$  if the cost functions are nonnegative.    If the jump Markov process, defined by the policy $\varphi$ and by the initial state distribution $\gamma,$ is nonexplosive, then the corresponding values of objective functions  for polices $\pi$ and $\varphi$ coincide without the assumption that the cost functions are nonnegative.  These facts hold for problems with multiple criteria.

If the transition rates $q(z,a)$ are bounded on $Gr(A)$, as this takes place in many applications to queueing control, it is well-known that the jump process under every policy is nonexplosive. However, jump processes under all policies may be nonexplosive even for problems with unbounded jump rates.  For example, jumps occur at arrival and departure epochs  in many controlled queues.  If the  rate of the arrival process, which may be Poisson or Markov-modulated, is bounded, then the total number of arrivals and departures over every finite deterministic interval of time   is finite with probability 1, and the corresponding Markov processes are nonexplosive even if their transition rates are unbounded. For example, transition rates can be unbounded because the number of servers is unbounded \cite{AKW,FZ} or because customers are impatient and can abandon the queue \cite{BSp}, but the corresponding jump Markov processes are nonexplosive.
 In general,  if the transition rates are unbounded, the corresponding jump process may be explosive.  \citet[Theorem 1]{PZ1} provided a general sufficient condition for the nonexplosiveness of  jumps processes defined by all policies.  
\subsection{Finite-horizon CTJMDPs.} 
For an initial state distribution $\gamma$ and a policy $\pi,$ the \textit{finite-horizon expected total discounted cost} with the discount rate $\alpha\in \R$ up to time $T \in \R_+$ is
\begin{equation}
\label{FDR}
V_{\alpha}^{T}(\gamma, \pi) :=  \BB{E}_\gamma^\pi \left[\int_{0}^{T \wedge t_\infty} e^{-\alpha s}c(\xi_{s}, \pi_s)ds + \sum\limits_{i=1}^{\infty} e^{-\alpha {u_i}}G_i(\xi_{u_i}, \pi_{u_i})\right] ,
\end{equation}
where the first summand is the total discounted cost collected up to the time $T$ with the cost rate $c(\cdot,\cdot)$ and  the second summand is the total discounted costs $G_i$ incurred at certain time epochs $(u_i\in [0,T])_{i=1,2,\ldots}.$   The functions $c:\bar{X}\times\bar{A}\mapsto \bar{\R}$ and $G_i:\bar{X}\times\bar{A}\mapsto \bar{\R}$ are assumed to be measurable with $c(x_\infty,\cdot)=c(\cdot,a_\infty)=G_i(x_\infty,\cdot)=G_i(\cdot,a_\infty)=0$ for all $i=1,2,\ldots.$ To avoid undefined sums, integrals, and expectations in \eqref{FDR}, we start with nonnegative functions $c$ and $G_i,$ $i=1,2,\ldots.$ We recall that, according to \eqref{Ext},  $f(\xi_{t}, \pi_t):=\int_{A(\xi_t)}f(\xi_{t}, a)\pi_t(da).$ The second summand in \eqref{FDR} models the situation when at some time instances $u_1,u_2,\ldots$ the decision maker has to make payments.  In particular, if $u_1=T$ and $G_i\equiv 0$ for $i=2,3,\ldots,$ then we have a problem with the terminal cost $G_1(\xi_T,\pi_T).$
\begin{theorem}
\label{thm:S1}
  Let the functions $c$ and $G_i,$ $i=1,2,\ldots,$ take nonnegative values, and let $T\in \R_+.$  For an initial distribution $\gamma$ on $X$ and a policy $\pi$, let $\varphi$ be a Markov policy satisfying \eqref{formula-de2018}. Then for all $\alpha \in \BB{R}$
\begin{equation}\label{eq61}
V_\alpha^T(\gamma, \varphi) \le V_\alpha^T(\gamma, \pi).\end{equation} 
If, in addition,  $P_\gamma^\varphi(T,X) =1,$   then $V_\alpha^T(\gamma, \varphi) = V_\alpha^T(\gamma, \pi).$
\end{theorem}
\proof{Proof.}
 Observe that
\begin{multline}
\label{eq:cost}
V_\alpha^T(\gamma, \pi)
= \int_0^{T} e^{-\alpha s}\BB{E}_\gamma^\pi c(\xi_{s}, \pi_s)ds  + \sum_{i=1}^\infty e^{-\alpha u_i}\BB{E}_\gamma^\pi G_i(\xi_{u_i},\pi_{u_i})\\
= \int_0^{T} e^{-\alpha s} \left(\int_X \int_{A(z)}  c(z,a) P_\gamma^\pi(s, dz, da) \right)ds + \sum_{i=1}^\infty e^{-\alpha u_i} \int_X \int_{A(z)}G_i(z,a) P_\gamma^\pi(u_i,dz,da),
\end{multline}
where the first equality follows from \eqref{FDR} since $\xi_s=x_\infty$ for $s\ge t_\infty$ and $c(x_\infty,\cdot)=0;$ the second equality follows from \eqref{f-Ex}  with $f(\xi_{t}, \pi_t) = c(\xi_{t}, \pi_t),$ $f(\xi_{t}, \pi_t) = G_i(\xi_{t}, \pi_t),$ and $Z=X.$  Then
\begin{multline}
\label{eq:sufT}
V_\alpha^T(\gamma, \varphi) = \int_0^{T} e^{-\alpha s} \int_X \int_{A(z)}  c(z,a) P_\gamma^\varphi(s, dz, da) ds + \sum_{i=1}^\infty e^{-\alpha u_i} \int_X \int_{A(z)}G_i(z,a) P_\gamma^\varphi(u_i,dz,da)  \\
\le \int_0^{T} e^{-\alpha s} \int_X \int_{A(z)}  c(z,a) P_\gamma^\pi(s, dz, da)ds + \sum_{i=1}^\infty e^{-\alpha u_i} \int_X \int_{A(z)}G_i(z,a) P_\gamma^\pi(u_i,dz,da)= V_\alpha^T(\gamma, \pi),
\end{multline}
where the first and last equalities follow from \eqref{eq:cost} applied to the policies $\varphi$ and $\pi$ respectively, and the inequality follows from  Theorem~\ref{main}. Thus \eqref{eq61} is proved.

Now let $P_\gamma^\varphi(T,X) = 1.$  Then    Theorem~\ref{main} implies that the inequality in \eqref{eq:sufT}  is an equality.
 \hfill\Halmos\endproof

Let $\Pi$ and $\Pi^M$ be respectively the classes of all history-dependent and Markov policies. The following corollary follows from Theorem~\ref{thm:S1}. 
\begin{corollary}
\label{C:suf}  Let the functions $c$ and $G_i,$ $i=1,2,\ldots,$ take nonnegative values, and let $T \in \R_+.$
 For an initial distribution $\gamma$ and
discount rate $\alpha \in \R$
\begin{equation}
\label{suf:inf}
\inf_{\varphi \in \Pi^M} V_\alpha^T(\gamma, \varphi) =  \inf_{\pi \in \Pi}V_\alpha^T(\gamma, \pi). 
\end{equation}
If, in addition,   $P_\gamma^\varphi(T,X) = 1$ for every Markov policy $\varphi,$ then
\begin{equation}
\label{suf:sup}
\sup_{\varphi \in \Pi^M} V_\alpha^T(\gamma, \varphi) =  \sup_{\pi \in \Pi}V_\alpha^T(\gamma, \pi). 
\end{equation}
\end{corollary}
\proof{Proof.}
 The inclusion $\Pi^M \subset \Pi$ implies $\inf_{\varphi \in \Pi^M} V_\alpha^T(\gamma, \varphi) \ge \inf_{\pi \in \Pi} V_\alpha^T(\gamma, \pi).$ The opposite inequality follows from Theorem~\ref{thm:S1}.
Equality \eqref{suf:sup} follows from the last statement of Theorem~\ref{thm:S1}.
 \hfill\Halmos\endproof

Now let us consider cost functions $c:\bar{X}\times\bar{A}\mapsto \bar{\R}$ and $G_i:\bar{X}\times\bar{A}\mapsto \bar{\R}$ without the assumption that they take nonnegative values.  For an arbitrary $r\in \bar{\R},$ let $r^+:=\max\{r,0\}$ and $r^-:=\min\{r,0\}$ be positive and negative parts of $r.$  Let $V_\alpha^{T,\oplus}(\gamma, \pi)$ and $V_\alpha^{T,\ominus}(\gamma, \pi)$ be the expected total costs defined in \eqref{FDR} with the cost functions $c$ and $G_i$ substituted with the functions $c^+,$ $G_i^+$ and $c^-,$ $G_i^-$ respectively, $i=1,2,\ldots.$  These  definitions imply  $V_\alpha^{T,\oplus}(\gamma, \pi)\ge 0$  and $V_\alpha^{T,\ominus}(\gamma, \pi)\le 0.$  If either $V_\alpha^{T,\oplus}(\gamma, \pi)<+\infty$ or $V_\alpha^{T,\ominus}(\gamma, \pi)>-\infty,$ then we say that the value $V_\alpha^T(\gamma, \pi)$ is \emph{defined} and set
\begin{equation}\label{eqsufdef}
V_\alpha^T(\gamma, \pi):=V_\alpha^{T,\oplus}(\gamma, \pi)+V_\alpha^{T,\ominus}(\gamma, \pi).
\end{equation}
For example, the values $V_\alpha^T(\gamma, \pi)$ are defined for all initial state distributions $\gamma,$ all policies $\pi,$ and all discount rates $\alpha\in\R$ if $T<+\infty,$  all the functions $c,$ $G_i$ are bounded either from below or from above simultaneously, and  for some natural number $k$ the functions $G_i$ are identically equal to 0 for $i>k.$  Of course, in this case \eqref{eqsufdef} holds.

The following two corollaries imply that, if the values $V_\alpha^T(\gamma, \pi)$ are defined, then under the assumption $P_\gamma^\varphi(T,X) = 1$ the corresponding conclusions of Theorem~\ref{thm:S1} and Corollary~\ref{C:suf} hold without the assumptions that the functions $c(x,a)$ and $G_i(x,a)$ take nonnegative values.
\begin{corollary}
\label{C:arbcostsfin} Let $T\in\R_+.$  For an initial distribution $\gamma$ on $X$ and a policy $\pi$, let $\varphi$ be a Markov policy satisfying \eqref{formula-de2018}. If the value $V_\alpha^T(\gamma, \pi)$ is defined and  $P_\gamma^\varphi(T,X) =1,$ then $V_\alpha^T(\gamma, \varphi) = V_\alpha^T(\gamma, \pi).$
\end{corollary}
\proof{Proof.}
The last claim of Theorem~\ref{thm:S1} applied to the functions $c^+,$ $G_i^+$ and $-c^-,$ $-G_i^-,$ $i=1,2,\ldots,$   implies respectively that $V_\alpha^{T,\oplus}(\gamma, \varphi)=V_\alpha^{T,\oplus}(\gamma, \pi)$ and $V_\alpha^{T,\ominus}(\gamma, \varphi)=V_\alpha^{T,\ominus}(\gamma, \pi).$ Thus,  the corollary follows from \eqref{eqsufdef}.
 \hfill\Halmos\endproof
 \begin{corollary}
\label{C:infsup}Let $T\in\R_+.$  For an initial distribution $\gamma$ on $X,$ if the value $V_\alpha^T(\gamma, \pi)$ is defined for every policy $\pi$ and  $P_\gamma^\varphi(T,X) =1$ for every Markov policy $\varphi,$ then equalities \eqref{suf:inf} and \eqref{suf:sup} hold.
\end{corollary}
\proof{Proof.}
This corollary follows from Corollary~\ref{C:arbcostsfin}.
\hfill\Halmos\endproof

Theorem~\ref{thm:S1} can be also applied to problems with multiple criteria and constraints.  Let, for a fixed initial state distribution $\gamma,$ the performance of policy $\pi$ is evaluated by a finite or infinite number $g(\gamma,\pi)\in\bar{\R}.$ In addition, for some nonempty set $\mathbb{B}$ there are a collection of functions $\{g_\beta:\P(X)\times\Pi\mapsto \bar{\R},\beta\in\mathbb{B}\}$ and a collection of real numbers $\{M_\beta:\beta\in\mathbb{B}\}.$ The constrained optimization problem is
\begin{equation}
\label{p:CP}
\begin{aligned}
{\rm minimize}_{\pi\in\Pi}\quad &g(\gamma, \pi)\\
\text{ subject to } \quad &g_\beta(\gamma, \pi) \le M_\beta, \qquad  \beta\in\mathbb{B}.
\end{aligned}
\end{equation}

Let us consider the finite time horizon $T\in \R_+,$ discount rate $\alpha\in\R,$ measurable cost functions $c:\bar{X}\times\bar{A}\mapsto \bar{\R}$ and $G_i:\bar{X}\times\bar{A}\mapsto \bar{\R},$ and time instances $u_i\in [0,T],$ $i=1,2,\ldots,$ satisfying the properties assumed in the first paragraph of this subsection.  For each $\beta\in\mathbb{B},$  let us consider the similar objects indexed by $\beta$ and satisfying the same properties. In particular, we consider the finite time horizons $T_\beta\in ]0,+\infty[,$ discount rates  $\alpha_\beta\in\R,$ measurable cost functions $c_\beta:\bar{X}\times\bar{A}\mapsto \bar{\R}$ and $G_{i,\beta}:\bar{X}\times\bar{A}\mapsto \bar{\R},$ and time instances $u_{i,\beta}\in [0,T_\beta],$ $i=1,2,\ldots.$ Let $V_{\alpha_\beta,\beta}^{T_\beta}(\gamma,\pi)$ denotes the expected total discounted costs defined in \eqref{FDR} for $\alpha:=\alpha_\beta,$ $c=c_\beta,$ $G_i:=G_{i,\beta},$ and  $u_i:=u_{i,\beta},$ where $\beta\in\mathbb{B},$ $i=1,2,\ldots.$ Let us consider the following assumption.
\begin{Assumption}
\label{Awd} Let the following conditions hold for a given initial distribution $\gamma$ on $X:$
\begin{enumerate}[{(i)}]
\item  either  the functions $c, G_i,$ $i=1,2,\ldots,$ take nonnegative values or the following two conditions hold: $P_\gamma^\varphi(T,X)=1$ for all Markov policies $\varphi,$ and  $V_\alpha^T(\gamma, \pi)$ is defined for all policies $\pi\in\Pi;$
\item  for each $\beta\in\mathbb{B}$ either the functions $c_\beta, G_{i,\beta},$ $i=1,2,\ldots,$ take nonnegative values or the following two conditions hold: $P_\gamma^\varphi(T_\beta,X)=1$ for all Markov policies $\varphi,$ and  $V_{\alpha_\beta}^{T_\beta}(\gamma, \pi)$ is defined for all policies $\pi\in\Pi.$
\end{enumerate}
\end{Assumption}

We remark that condition (ii) in Assumption~\ref{Awd} is condition (i) applied to  a finite-horizon CTJMDP with $c,$ $G_i,$ $\alpha,$ and $T$ replaced with $c_\beta,$ $G_{i,\beta},$ $\alpha_\beta,$ and $T_\beta$ respectively.
\begin{corollary}\label{constrprobfg}
For an  initial state distribution $\gamma$ on $X,$ let us consider problem \eqref{p:CP} with $g(\gamma,\pi)=V_\alpha^T(\gamma, \pi)$ and $g_\beta(\gamma,\pi)=V_{\alpha_\beta,\beta}^{T_\beta}(\gamma,\pi)$ for all $\pi\in\Pi$ and $\beta\in\mathbb{B}.$  If Assumption~\ref{Awd} holds for the initial distribution $\gamma,$ then for every feasible policy $\pi$ inequality \eqref{eq61} holds for a Markov policy $\varphi$ satisfying \eqref{formula-de2018}, and the policy $\varphi$ is feasible.
\end{corollary}
\proof{Proof.}
Let us consider problem \eqref{p:CP} with $g(\gamma,\pi)=V_\alpha^T(\gamma,\pi)$ and $g_\beta(\gamma,\pi)=V_{\alpha_\beta,\beta}^{T_\beta}(\gamma,\pi).$ For an arbitrary feasible policy $\pi,$ let us consider a Markov policy $\varphi$ satisfying \eqref{formula-de2018}. Theorem~\ref{main} and Corollary \ref{C:arbcostsfin} imply that
$V_\alpha^T(\gamma,\varphi)\le V_\alpha^T(\gamma,\pi)$ and $V_{\alpha_\beta,\beta}^{T_\beta}(\gamma,\varphi)\le V_{\alpha_\beta,\beta}^{T_\beta}(\gamma,\pi)\le M_\beta$  for  all $\beta\in\mathbb{B}.$
\hfill\Halmos\endproof

  Thus, for every feasible policy, that is, a policy satisfying constraints in \eqref{p:CP}, there is a feasible Markov policy with the same or smaller objective function.  The same is true  if the objective function and the functions in constraints are sums of finite numbers of the expected discounted costs with possibly different discount rates and horizons and  if Assumption~\ref{Awd} is satisfied for all the expected discounted total costs in the sums.  The latter means that the summands in the objective function satisfy Assumption~\ref{Awd}(i) and the summands in constraints satisfy Assumption~\ref{Awd}(ii).

\subsection{Infinite horizon  CTJMDPs.}  For infinite-horizon problems, in addition to cost rates $c$ and instant costs $G_i,$ we shall also consider costs $C(\xi_{t_{n-1}}, \xi_{t_n})$ incurred at jump epochs $t_n,$ $n=1,2,\ldots.$
The cost structure of an infinite-horizon CTJMDP is defined by the following three nonnegative cost functions:
\begin{itemize}
\item[(i)] the cost rate function $c(z,a)$ representing the cost per unit time when an action $a$ is chosen at state $z$;

\item[(ii)] the instant cost function $G_i(z)$   representing  the costs collected at time epochs  $u_i,$ where $u_i\in\R_+^0,$ if $z=\xi_{u_i},$ $i=1,2,\ldots;$

\item[(iii)] the instantaneous cost function $C(z,y)$   representing the cost incurred at the jump epoch when the process transitions from state $z$ to state $y$.
\end{itemize}

 The functions $c:\bar{X}\times\bar{A}\mapsto \bar{\R},$  $G_i:\bar{X}\times\bar{A}\mapsto \bar{\R},$ and $C:\bar{X}\times \bar{X}\mapsto \bar{\R}$ are assumed to be measurable with $c(x_\infty,\cdot)=c(\cdot,a_\infty)=G_i(x_\infty,\cdot)=G_i(\cdot,a_\infty)=C(x_\infty,\cdot)=C(\cdot,x_\infty)=0$ for all $i=1,2,\ldots.$
If all the cost functions $c,$ $G_i,$ $i=1,2,\ldots,$ and $C$ take either always nonnegative values or always nonpositive values,
%
for an initial state distribution $\gamma$, a policy $\bar{\pi}$, and a discount rate $\alpha \in \R_+^0,$ the \textit{infinite-horizon expected total discounted cost} is
\begin{equation}
\label{IDR}
V_{\alpha}(\gamma, \bar{\pi}): =  \BB{E}_\gamma^{\bar{\pi}} \left[\int_{0}^{t_\infty} e^{-\alpha s}c(\xi_{s}, \pi_s)ds + \sum\limits_{n=1}^\infty e^{-\alpha t_n}C(\xi_{t_{n-1}}, \xi_{t_n}) + \sum\limits_{i=1}^{\infty} e^{-\alpha {u_i}}G_i(\xi_{u_i}, \bar{\pi}_{u_i})\right].
\end{equation}
\begin{theorem}
\label{thm:S1inf}
  Let the functions $c,$ $C,$ and $G_i,$ $i=1,2,\ldots,$ take nonnegative values.  For an initial distribution $\gamma$ on $X$ and a policy $\pi$, let $\varphi$ be a Markov policy satisfying \eqref{formula-de2018}. Then for all $\alpha \in \BB{R}_+^0$
\begin{equation}\label{eq69}
V_\alpha(\gamma, \varphi) \le V_\alpha(\gamma, \pi).\end{equation} 
If, in addition,  $P_\gamma^\varphi(t,X) =1$ for all $t>0,$   then $V_\alpha(\gamma, \varphi) = V_\alpha(\gamma, \pi).$
\end{theorem}
\proof{Proof.} First, we prove the theorem for problems without instantaneous costs at jump epochs, that is $C\equiv 0.$  In this case, let us denote by $\bar{V}_\alpha(\gamma,\bar{\pi}),$ where $\bar{\pi}$ is an arbitrary policy, the expected total infinite-horizon cost defined in \eqref{IDR} with the omitted second summand in the right-hand side.   The  finite-horizon version of these costs up to the epoch $T\in \R_+$ is
\begin{equation}\label{eq70EF}
\bar{V}_\alpha^T(\gamma,\bar{\pi}):=\BB{E}_\gamma^{\bar{\pi}} \left[\int_{0}^{T\wedge t_\infty} e^{-\alpha s}c(\xi_{s}, \bar{\pi}_s)ds  + \sum\limits_{i=1}^{\infty} e^{-\alpha {u_i}}1\{u_i\le T\}G_i(\xi_{u_i}, \bar{\pi}_{u_i})\right].
\end{equation}
The monotone convergence theorem implies that $\bar{V}_\alpha^T(\gamma,\bar{\pi})\uparrow \bar{V}_\alpha(\gamma,\bar{\pi})$ as $T\to +\infty.$ Therefore,
\begin{equation}\label{eineqefefFeller}
\bar{V}_\alpha(\gamma,\varphi)=\lim_{T\to\infty} \bar{V}_\alpha^T(\gamma,\varphi)\le \lim_{T\to\infty} \bar{V}_\alpha^T(\gamma,\pi)=\bar{V}_\alpha(\gamma,\pi),
\end{equation}
where the inequality follows from Theorem~\ref{thm:S1}. In addition, as follows from Theorem~\ref{thm:S1}, if $P_\gamma^\varphi(T,X) =1$ for all $T>0,$ then the inequality in
\eqref{eineqefefFeller} holds in the form of an equality.    The theorem is proved for $C\equiv 0.$  To complete the proof of the theorem, it is sufficient to show that
\begin{equation}\label{e71}
\BB{E}_\gamma^\varphi \sum\limits_{n=1}^\infty e^{-\alpha t_n}C(\xi_{t_{n-1}}, \xi_{t_n})\le \BB{E}_\gamma^\pi \sum\limits_{n=1}^\infty e^{-\alpha t_n}C(\xi_{t_{n-1}}, \xi_{t_n})
\end{equation}
and,  if $P_\gamma^\varphi(T,X) =1$ for all $T>0,$ then the inequality in
\eqref{e71} holds in the form of an equality.

To prove \eqref{e71}, we set $\tilde{C}(x_\infty,\delta_{a_\infty}):=0$ and consider the function $\tilde{C}:X\times\P(A)\mapsto [0,+\infty],$
\begin{equation}\label{eq73}\tilde{C}(z,p):=\int_{X\setminus \{z\}}C(z,y)\tilde{q}(z,p,dy),\end{equation}
where the measures $\tilde{q}(z,p,dy)$ and $\tilde{q}(z,p)$ are defined in  \eqref{q-def21} and \eqref{q-def22} respectively. Let us fix $n=1,2,\ldots$ and consider the nonnegative random variable $s_n:=t_n-t_{n-1},$ if $t_{n-1}<+\infty,$ which is the sojourn time. We also set $s_n:=+\infty$ if $t_{n-1}=+\infty.$  As follows from \eqref{JMP}, \eqref{policy}, and \eqref{m2}, for $t\in \bar{\R}_+^0,$
\begin{equation}\label{eq76}
\BB{P}_\gamma^{\bar{\pi}}\{s_n\le t|\,\cal{F}_{t_{n-1}}\}=1-\exp(-\int_0^t\tilde{q}(\xi_{t_{n-1}},\bar{\pi}^{n-1}(x_0,t_1,x_1,\ldots,t_{n-1},x_{n-1},s))ds),
\end{equation}
where the function $\tilde{q}:X\times\P(A)\mapsto [0,+\infty]$ is defined in \eqref{q-def22} and $\tilde{q}(x_\infty,\delta_{a_\infty}):=0.$ Since the state $x_\infty$ is always absorbing,  we also  set $\tilde{q}(x_\infty,\delta_{a_\infty},B):=0$ for $B\in \B(\bar{X}).$  As follows from \eqref{JMP}, \eqref{policy}, \eqref{m2}, and \eqref{eq76}, for $B\in\B(X),$ we have that $\BB{P}_\gamma^{\bar{\pi}}\{s_n< +\infty,\,\tilde{q}(\xi_{t_{n-1}},\bar{\pi}_{t_{n-1}+s_n}) =0 |\,\cal{F}_{t_{n-1}}\}=0$ and, following everywhere the convention $\frac{0}{0}:=0,$
\begin{equation}\label{eq75EF}
\BB{P}_\gamma^{\bar{\pi}}\{ \xi_{t_n}\in B|{\cal F}_{t_{n-1}}\}=\BB{E}_\gamma^{\bar{\pi}}\left[\frac{\tilde{q}(\xi_{t_{n-1}},\bar{\pi}_{t_{n-1}+s_n},B\setminus\{\xi_{t_{n-1}}\})}{\tilde{q}(\xi_{t_{n-1}},\bar{\pi}_{t_{n-1}+s_n})}
{\big|}\,\cal{F}_{t_{n-1}}\right].
\end{equation}
Then
\begin{equation}\label{eq75}
\begin{aligned}
&\BB{E}_\gamma^{\bar{\pi}} [e^{-\alpha t_n}C(\xi_{t_{n-1}},\xi_{t_n})|\,\cal{F}_{t_{n-1}}]=
e^{-\alpha t_{n-1}}\BB{E}_\gamma^{\bar{\pi}} \left[e^{-\alpha s_n}\frac{\tilde{C}(\xi_{t_{n-1}},\bar{\pi}_{t_{n-1}+s_n})}{\tilde{q}(\xi_{t_{n-1}},\bar{\pi}_{t_{n-1}+s_n})}{\big|}\, \cal{F}_{t_{n-1}}\right]\\
&=e^{-\alpha t_{n-1}}\int_0^\infty e^{-\alpha s}\frac{\tilde{C}(\xi_{t_{n-1}},\bar{\pi}_{t_{n-1}+s})}{\tilde{q}(\xi_{t_{n-1}},\bar{\pi}_{t_{n-1}+s})}d\BB{P}_x^{\bar{\pi}}\{s_n\le s|\,\cal{F}_{t_{n-1}}\}\\
&=e^{-\alpha t_{n-1}}\int_0^\infty e^{-\alpha s}\tilde{C}(\xi_{t_{n-1}},\bar{\pi}_{t_{n-1}+s})\BB{P}_x^{\bar{\pi}}\{s_n> s|\,\cal{F}_{t_{n-1}}\}dt,
\end{aligned}
\end{equation}
where the first equality in \eqref{eq75} follows from \eqref{eq75EF} and because the random variable $t_{n-1}$ is $\cal{F}_{t_{n-1}}$-measurable, the second equality holds because conditional expectation can be written as an integral with respect to the conditional distribution, and the last equality follows from \eqref{eq76} and from the explicit differentiation in $s$ the function $\BB{P}_x^{\bar{\pi}}\{s_n\le s|\,\cal{F}_{t_{n-1}}\}$ displayed in \eqref{eq76} with $s=t.$

Let us consider the nonnegative function $f(s)=e^{-\alpha s}\tilde{C}(\xi_{t_{n-1}},\bar{\pi}_{t_{n-1}+s}).$ According to \cite[p. 263]{Fei94},
$\BB{E}_\gamma^{\bar{\pi}}[\int_0^{s_n} f(s)ds|\,\cal{F}_{t_{n-1}}]=\int_0^\infty f(s)\BB{P}_x^{\bar{\pi}}\{s_n> s|\,\cal{F}_{t_{n-1}}\}ds.$  This formula and \eqref{eq75} imply
\begin{equation}\label{eq755}
\BB{E}_\gamma^{\bar{\pi}} [e^{-\alpha t_n}C(\xi_{t_{n-1}},\xi_{t_n})|\,\cal{F}_{t_{n-1}}]=\BB{E}_\gamma^{\bar{\pi}}[\int_0^{s_n} e^{-\alpha (t_{n-1}+ s)}\tilde{C}(\xi_{t_{n-1}},\bar{\pi}_{t_{n-1}+s})ds|\,\cal{F}_{t_{n-1}}].
\end{equation}
   By changing the variable $s$ to $t:=t_{n-1}+s$ in \eqref{eq755} and taking expectations in  \eqref{eq755}, we have

\begin{equation}\label{eq72}
\BB{E}_\gamma^{\bar{\pi}} e^{-\alpha t_n}C(\xi_{t_{n-1}}, \xi_{t_n})=\BB{E}_\gamma^{\bar{\pi}} \int_{t_{n-1}}^{t_n} e^{-\alpha t} \tilde{C}(\xi_{t_{n-1}},\bar{\pi}_t)dt, \qquad n=1,2,\ldots,\ \bar{\pi}\in\Pi,
\end{equation}
which implies
\[
\BB{E}_\gamma^{\bar{\pi}}\sum_{n=1}^\infty e^{-\alpha t_n}C(\xi_{t_{n-1}}, \xi_{t_n})=\BB{E}_\gamma^{\bar{\pi}}\int_0^{t_\infty} e^{-\alpha t} \tilde{C}(\xi_{t},\bar{\pi}_t)dt, \qquad \bar{\pi}\in\Pi.
\]
Thus, the expected discounted sum of instant costs $C(z,y)$ at jump epoch is equal to the expected total discounted cost  with the cost rate $\tilde{C}(z,a)=\int_{X\setminus\{z\}}C(z,y){\tilde q}(z,a,dy).$  Since the theorem is proved in  \eqref{eineqefefFeller}  for problems without instant costs at jump epochs,  equality \eqref{eq72} implies inequality \eqref{e71}.  In addition, as follows from Theorem~\ref{thm:S1}, if $P_\gamma^\varphi(T,X) =1$ for all $T>0,$ then the inequality in
\eqref{e71} holds in the form of an equality.
\hfill\Halmos\endproof

\begin{corollary}
\label{C:sufinf}  Let the functions $c,$ $C,$ and $G_i,$ $i=1,2,\ldots,$ take nonnegative values.
 For an initial distribution $\gamma$ and
discount rate $\alpha \in \R_+^0$
\begin{equation}
\label{suf:infinf}
\inf_{\varphi \in \Pi^M} V_\alpha(\gamma, \varphi) =  \inf_{\pi \in \Pi}V_\alpha(\gamma, \pi). 
\end{equation}
If, in addition,   $P_\gamma^\varphi(t,X) = 1$ for every Markov policy $\varphi$ and for every $t>0,$ then
\begin{equation}
\label{suf:supinf}
\sup_{\varphi \in \Pi^M} V_\alpha(\gamma, \varphi) =  \sup_{\pi \in \Pi}V_\alpha(\gamma, \pi). 
\end{equation}
\end{corollary}
\proof{Proof.}  The proof is identical to the proof of Corollary~\ref{C:suf} with Theorem~\ref{thm:S1inf} used instead of Theorem~\ref{thm:S1}. \hfill\Halmos\endproof

If the cost functions $c,$ $G_i,$  and $C$ can take positive and negative values, we consider the infinite-horizon expected total discounted costs
$V_{\alpha}^\oplus(\gamma, \pi)$ and $V_{\alpha}^\ominus(\gamma, \pi)$ for costs functions $c^+,$ $G^+_i,$  $C^+$ and $c^-,$ $G^-_i,$  and $C^-$ respectively, $i=1,2,\ldots.$
We say that the infinite-horizon expected total discounted cost $V_{\alpha}(\gamma, \pi)$ is defined if either $V_{\alpha}^\oplus(\gamma, \pi)<+\infty$ or $V_{\alpha}^\ominus(\gamma, \pi)>-\infty.$ If $V_{\alpha}(\gamma, \pi)$ is defined, we set similarly to \eqref{eqsufdef}
\begin{equation}\label{eqsufdefinf}
V_\alpha(\gamma, \pi):=V_\alpha^\oplus(\gamma, \pi)+V_\alpha^\ominus(\gamma, \pi).
\end{equation}
The following two corollaries are the infinite-horizon versions of Corollaries~\ref{C:arbcostsfin} and \ref{C:infsup}, and they follow from Theorem~\ref{thm:S1inf} in the same way as Corollaries~\ref{C:arbcostsfin} and \ref{C:infsup} follow from Theorem~\ref{thm:S1}.
\begin{corollary}
\label{C:arbcostsinfininf}  For an initial distribution $\gamma$ on $X$ and a policy $\pi$, let $\varphi$ be a Markov policy satisfying \eqref{formula-de2018}. If the value $V_\alpha(\gamma, \pi)$ is defined and  $P_\gamma^\varphi(t,X) =1$  for all $t>0, $ then $V_\alpha(\gamma, \varphi) = V_\alpha(\gamma, \pi).$
\end{corollary}
 \begin{corollary}
\label{C:infsupinf} For an initial distribution $\gamma$ on $X,$ if the value $V_\alpha(\gamma, \pi)$ is defined for every policy $\pi$ and  $P_\gamma^\varphi(t,X) =1$ for every $t>0$ and for every Markov policy $\varphi,$ then equalities \eqref{suf:infinf} and \eqref{suf:supinf} hold.
\end{corollary}

For an infinite horizon,  \textit{the average cost per unit time } is
\begin{equation}
\label{AVGR}
W(\gamma,\pi) := \limsup_{\alpha \downarrow 0}\alpha V_{\alpha}(\gamma,\pi). 
\end{equation}
Another way  to define the average cost per unit time is
\begin{equation}
\label{AVGR1}
W^1(\gamma, \pi):= \limsup_{T \to +\infty} \frac{V_{0}^{T}(\gamma, \pi)}{T},
\end{equation}
where   $V_{0}^{T}(\gamma, \pi)$ is introduced in \eqref{IDR}.  The average cost $W(\gamma,\pi)$ is defined if  $V_{\alpha}(\gamma,\pi)$ is defined for each $\alpha >0.$  The average cost  $W^1(\gamma,\pi)$ is defined if $\bar{V}_{0}^{T}(\gamma, \pi)$ is defined for each $T>0.$  These definitions are natural  if the initial distribution $\gamma$ and the policy $\pi$ define an nonexplosive process, that is, $P_\gamma^\pi (t,X)=1$ for all $t>0.$

Average costs \eqref{AVGR} and \eqref{AVGR1} are related under certain conditions. If instant and jump costs are equal to zero and the stochastic process $\xi_t$ defined by an initial state distribution $\gamma$ and a policy $\pi$  is nonexplosive, that is, $C=G_i\equiv 0,$ $i=1,2,\ldots,$ and $P_\gamma^\pi\{\xi_t\in X\}=1$ for all $t>0,$ then
\begin{equation}
\label{IDR0}
V_{\alpha}(\gamma, \pi):=  \BB{E}_\gamma^\pi \int_{0}^{+\infty} e^{-\alpha s}c(\xi_{s}, \pi_s)ds=\int_{0}^{+\infty} e^{-\alpha s}\BB{E}_\gamma^\pi[c(\xi_{s}, \pi_s)]ds.
\end{equation}
%
%
In this case, as follows from the Tauberian  theorem \cite[p. 197]{GHL}, $W(\gamma, \pi) \le W^1(\gamma, \pi)$ if the function $c$ takes nonnegative values. In addition, this inequality holds as an equality  if either of the limits $\lim_{T \to \infty} \frac{V_{0}^{T}(\gamma, \pi)}{T}$
or $\lim_{\alpha \downarrow 0} \alpha V_{\alpha}(\gamma,\pi)$ exists and is finite; \cite[Chapter V, Corollary 1a and Theorem 14]{Widder}.

 Standard properties of limits, Theorem~\ref{thm:S1inf}, and Corollaries~\ref{C:sufinf}--\ref{C:infsupinf} imply that the statements of Theorem~\ref{thm:S1inf} and Corollaries~\ref{C:sufinf}--\ref{C:infsupinf} remain valid, if average costs per time $W$ and $W^1$ are used instead of the expected discounted costs $V_\alpha$ and the corresponding objective criteria are defined.  In addition, this is true for the finite sums of the objective criteria $V_\alpha(\gamma,\pi),$ $W(\gamma,\pi),$ and $W^1(\gamma,\pi)$  with possibly different cost functions and for possibly different discount factors for the total discounted criteria if each summand in each sum is defined.  In particular, if, for a given initial distribution $\gamma,$  every Markov policy defines a nonexplosive jump Markov process, then for every feasible policy for problem \eqref{p:CP} with $g$ and $g_\beta$ being such finite sums, there exists a feasible Markov policy with the equal or smaller objective function.

 The CTJMDP literature usually deals only with cost rate functions $c$ and, if the problem is finite-horizon, with  terminal costs. The instantaneous costs $G_i$ has never been considered before. Piunovsky and Zhang \cite{P1} considered instantaneous costs $C(z)$ at jump epochs that depend only on the state $z$ from which the process jumps. Feinberg~\cite{Fei04, Fei12} considered instantaneous costs at jump epochs of the form $C(z,a,y),$ where $z$ is the state from which the process jumps, $a$ is the action selected at state $z$ at the jump epoch, and $y$ is the state to which the process jumps. It is shown in \cite[Corollary 4.4]{Fei04} that for nonrelaxed (nonrandomized) policies the instantaneous jump costs $C(z,a,y)$ can be substituted with the  additional cost rate
 \begin{equation}\label{etrans} \bar{C}(z,a,y):=\int_X C(z,a,y)\tilde{q}(z,a,dy)
  \end{equation}
  added to the cost rate $c(z,a).$ 
 The results for general (relaxed) policies stated in \cite{Fei04, Fei12} have correct proofs only for the instantaneous costs of the form $C(z,y),$ that is, for instantaneous costs at jump epochs that do not depend on actions.  

 In particular, it was overlooked in \cite[p. 510]{Fei04} that for general (relaxed) policies the  transformation in \eqref{etrans} leads to the instant cost  $\tilde{C}(z,p,y):=\int_{A(z)}\int_X C(z,a,y)\bar{q}(z,a,dy)p(da),$ where $p\in\mathcal{P}(A(z))$ is a relaxed action, rather than to the desired costs $\bar{C}(z,p,y):= \int_X C(z,p,y)\tilde{q}(z,p,dy)$ with $C(z,p,y):=\int_{A(z)} C(z,a,y)p(da)$  and with $\tilde{q}(z,p,dy)$ defined in \eqref{q-def21}.  However, for $C(z,a,y):=C(z,y),$ the equality $\int_{A(z)}\int_X C(z,y)\tilde{q}(z,a,dy)p(da)= \int_X C(z,y)\tilde{q}(z,p,dy)$ holds.  That is, the desired transformation takes place in \cite{Fei04, Fei12} if the values $C(z,a,y)$ do not depend on $a.$

In the following example instantaneous costs $C(z,a,y)$ depend on the action $a,$ all Markov policies define nonexplosive Markov chains, and the conclusions of Theorems~\ref{thm:S1} and \ref{thm:S1inf} do not hold.

\begin{example}
\label{Ex-C}
{\rm {\it  For an initial distribution $\gamma$ and for a policy $\pi$, the Markov policy $\varphi$ satisfying \eqref{formula-de2018} can have higher infinite-horizon expected total discounted costs than the policy $\pi$ when the instantaneous costs $C$ depend on the action chosen at the jump epoch.}
Let $X = \{1,2\},$ $ A = \{ b,c\},$ $ A(1) = b,$ $ A(2) = A,$ $ \tilde{q}(1,b) = \tilde{q}(2,b) = 2$, and $\tilde{q}(2,c) = 1$. The cost rate function $c(x,a) = 0$ for all $x \in X$ and $a \in A$, and the instantaneous costs $C$ are $C(1,b,2) = 0,$ $ C(2,b, 1) = 1,$ and $C(2,c,1) = 2$. There are no instant costs $G_i,$ that is, $G_i\equiv 0$ for all $i=1,2,\ldots.$ This CTJMDP is described in Figure~\ref{figure-2}.
\begin{figure}[h]
\centering
\begin{tikzpicture}[>=stealth',shorten >=1pt,auto,node distance=3.5cm]
  \node[state] (q1)      {$1$};
  \node[state]         (q2) [right of=q1]  {$2$};

  \path[->]          (q1)  edge   [bend left=70]  node[above,text width=2.7cm,align=center ]{$a=b$,$\tilde{q}(1,b) = 2,$ $C(1,b,2) = 0$} (q2); (q2);
  \path[->]          (q2)  edge   [bend left=70]   node[below,text width=2.7cm,align=center ]{$a=b$,$\tilde{q}(2,b) = 2,$ $C(2,b,1) = 1$} (q1);
  \path[->]          (q2) edge     node[above,align=center ] {$a=c$,$\tilde{q}(2,c) = 1,$} (q1);
 \path[->]          (q2) edge     node[below, align=center ] { $C(2,c,1) = 2$} (q1);
\end{tikzpicture}
\caption{CTJMDP with two states and two actions.}
\label{figure-2}
\end{figure}
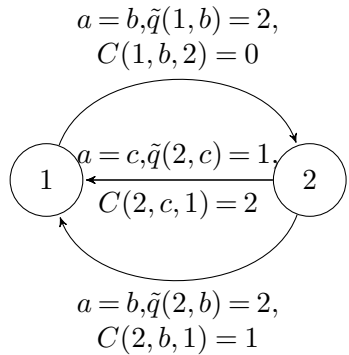

 Let $N(t,2)$ represent the number of jumps into state 2 up to time $t$ and let $\pi$ be a non-randomized policy choosing the action $\pi_t$ at time $t,$ where
\begin{equation}
\label{pi}
\pi_t = bI\{\xi_t = 1\} + b I\{\xi_t = 2, N(t,2) \text{ is even or } 0\} + cI\{\xi_t = 2, N(t,2) \text{ is odd}\}.
\end{equation}
Let the initial state be $2$. Observe that, $\BB{P}_2^\sigma(\xi_t \in X) = 1$ for all $t \in \BB{R}_+$ for every policy $\sigma$. This follows from Corollary~\ref{C:main}. Since $C(1,b,2)=0, $ the expected discounted total cost up to time $t$ under every policy  does not change in $t$ when the process is at state 1.  Since $C(2,b,1)\tilde{q}(2,b,1)=C(2,c,1)\tilde{q}(2,c,1)=2,$ then $C(2,a_t,1)\tilde{q}(2,a_t,1)=2$ for any nonrandomized policy. Therefore,
 the discounted total cost rate increases with the rate $2e^{-\alpha t}$ if the policy $\pi$ is used at state 2 at time $t.$ Thus,
 \begin{equation}\label{eq86} V_\alpha(2,\pi)=\BB{E}_2^\pi\int_0^{+\infty} 2e^{-\alpha t}I\{\xi_t=2\}dt=2\int_0^{+\infty} e^{-\alpha t}P_2^\pi(t,2)dt.\end{equation}

Let $\varphi$ be a Markov policy satisfying \eqref{formula-d} for all $t \in \BB{R}_+$. That is, the Markov policy $\varphi$ selects the action $b$ in state $1$ and an action $a \in \{b,c\}$ in state $2$ with probability $\frac{P_2^\pi(t, 2,a)}{P_2^\pi(t, 2)}$. Observe that $P_2^\pi(t, 2,a)>0$ and $\varphi_t(a|2,t)>0$ for $t>0,$ $a\in A(2)=\{b,c\}.$ For the Markov policy $\varphi,$ if the process is at state $2$ at time $t>0,$ then the jump rate is
\begin{equation}
\label{C-61}
\tilde{q}(2,\varphi_t,1)=\tilde{q}(2,b)\varphi(b|2,t) + \tilde{q}(2,c)\varphi(c|2,t)
=2\varphi(b|2,t) + \varphi(c|2,t) = 1+ \varphi(b|2,t),
\end{equation}
and  the expected instantaneous cost incurred, if a jump occurs at the epoch $t,$
is
\begin{equation}
\label{C-51}
C(2,\varphi_t,1) =C(2,b,1)\varphi(b|2,t) + C(2,c,1)\varphi(c|2,t)
=\varphi(b|2,t) + 2\varphi(c|2,t) = 1+\varphi(c|2,t).
\end{equation}
In view of \eqref{C-61} and  \eqref{C-51}, the expected discounted total cost  increases with the rate $e^{-\alpha t}(1+ \varphi(b|2,t))(1+\varphi(c|2,t))=e^{-\alpha t}(2+ \varphi(b|2,t)\varphi(c|2,t))$ if the policy $\varphi$ is used at state 2 at time $t.$

Therefore, starting from initial state 2, the infinite-horizon expected total discounted cost earned by the Markov policy $\varphi$  is
\begin{multline}
\label{V-varphi}
\begin{aligned}[t]
V_\alpha(2,\varphi) &=  \BB{E}_2^\varphi \int_{0}^{\infty} e^{-\alpha t}(2+\varphi(b|2,t)\varphi(c|2,t))I\{\xi_t = 2\} )dt \\
&= \int_0^{\infty} e^{-\alpha t}(2+ \varphi(b|2,t)\varphi(c|2,t)) P_2^\varphi(t, 2)dt\\
&= 2\int_0^{\infty} e^{-\alpha t}P_2^\pi(t, 2)dt + \int_0^{\infty} e^{-\alpha t} \varphi(b|2,t)\varphi(c|2,t) P_2^\pi(t, 2)dt > V_\alpha(2,\pi),
\end{aligned}
\end{multline}
where the inequality follows from \eqref{eq86} and  $ \int_0^{\infty} e^{-\alpha t} \varphi(c|2,t)\varphi(b|2,t) P_2^\pi(t, 2)dt>0$. Since in this example $V^T_\alpha(x,\sigma)\uparrow V_\alpha(x,\sigma)$ as $T\to\ +\infty$ for any policy $\sigma$ and for any initial state $x,$ then $V^T_\alpha(2,\varphi)>V^T_\alpha(2,\pi)$ in this example for sufficiently large $T.$
}
\end{example}

%
%

\section*{Acknowledgments.} Research of the first author was partially supported by the
National Science Foundation grant CMMI-1636193. Research of the third author was supported by the Russian Science
Foundation project 19-11-00290.
The authors thank Peng Dai, Pavlo O. Kasyanov, and Yi Zhang for valuable remarks.



\bibliographystyle{informs2014} 
\bibliography{references} 


\end{document}